\documentclass[a4paper,oneside,11pt]{article}

\usepackage[T1]{fontenc}             
\usepackage[french,english]{babel}           %
\usepackage[utf8]{inputenc}        
\usepackage{amsmath,amssymb,dsfont}  
\usepackage{textcomp}                
\usepackage{graphicx,color}          
\usepackage{fullpage}                
\usepackage{a4wide}                  
\usepackage[colorlinks,linkcolor=blue,citecolor=blue]{hyperref} 
\usepackage{ulem}                    
\usepackage{authblk}

\usepackage{amsthm}
\newtheorem{theorem}{Theorem}[section]
\newtheorem{proposition}{Proposition}[section]
\newtheorem{lemma}{Lemma}[section]
\newtheorem{corollary}{Corollary}[section]

\theoremstyle{definition}

\theoremstyle{remark}

\newcommand{\pa}[1]{\ensuremath{\left( #1 \right)}}
\newcommand{\cro}[1]{\ensuremath{\left[ #1 \right]}}
\newcommand{\ac}[1]{\ensuremath{\left\{ #1 \right\}}}
\newcommand{\abs}[1]{\ensuremath{\left| #1 \right|}}

\newcommand{\argmin}{\operatornamewithlimits{argmin}}

\newcommand{\KL}{\operatorname{KL}}
\newcommand{\TV}{\operatorname{TV}}

\newcommand{\N}{\ensuremath{\mathds{N}}}

\newcommand{\R}{\ensuremath{\mathds{R}}}

\renewcommand{\L}{\ensuremath{\mathds{L}}}

\newcommand{\proba}[1]{\ensuremath{\mathds{P}\!\left(#1\right)}}			
\newcommand{\esp}[1]{\ensuremath{\mathds{E}\left[  #1 \right]}} 		
\newcommand{\esps}[2]{\ensuremath{\mathds{E}_{#1}\!\left[ #2 \right]}}	

\newcommand{\var}[1]{\ensuremath{\operatorname{Var}\!\left(  #1 \right)}} 		

\newcommand{\1}[1]{\ensuremath{\mathds{1}_{ #1}}}		

\newcommand{\norm}[1]{\left\| #1 \right\|}

\newcommand{\mi}[1]{\boldsymbol{\underline{#1}}} 
\newcommand{\ml}{\mi{\ell}}
\newcommand{\mj}{\mi{j}}
\newcommand{\mcL}{\mathcal{L}}
\newcommand{\mcJ}{\mathcal{J}}
\renewcommand{\L}{\mathds{L}}

\begin{document}

\title{Minimax density estimation in the adversarial framework under local differential privacy}
\author[1]{M\'elisande Albert} 
\author[1]{Juliette Chevallier}
\author[1]{B\'eatrice Laurent}
\author[2]{Ousmane Sacko}
\affil[1]{Institut de Math\'ematiques de Toulouse; UMR 5219, Universit\'e de Toulouse, INSA}
\affil[2]{MODAL'X; UMR 9023, Université Paris Nanterre}
\maketitle

\renewcommand{\abstractname}{\vspace{-1cm}}
\begin{abstract}
We consider the problem of nonparametric density estimation under privacy constraints in an adversarial framework. To this end, we study minimax rates over Sobolev spaces under local differential privacy. We first obtain a lower bound  which allows us to quantify the impact of privacy compared with the classical framework. Next, we introduce a new Coordinate block privacy mechanism that guarantees local differential privacy, which, coupled with a projection estimator, achieves the minimax optimal rates. Finally, we develop an adaptive procedure which is optimal in the minimax sense up to logarithmic terms.
\end{abstract}

\bigskip


\medskip



\section{Introduction}
\label{sect:Intro}

In this paper, we study minimax nonparametric density estimation under local differential privacy constraints, in an adversarial framework. Let us first present each context.

\paragraph{Minimax nonparametric density estimation.}
Let $X_1, \dots, X_n$ be $n$  independent and identically distributed (i.i.d.) random variables with common density $f$ on $[0,1]^d$ with respect to the Lebesgue measure.  
We aim at estimating the underlying density $f$, a classical challenge in statistics extensively studied in the literature.
Early attempts to address this problem rely on kernel methods, as done by \cite{rosenblatt1956remarks,parzen1962estimation,silverman1978weak}. Projection methods also emerge as a viable approach, built upon the decomposition of $f$ into orthonormal bases. Notable contributions include the works of \cite{schwartz1967estimation,kronmal1968estimation,walter1977properties,efromovich1999nonparametric}.
We want to construct adaptive estimators in the minimax sense.
More precisely, consider a loss function $\ell$ and a regularity space $\mathcal{F}$. 
We say that an estimator $\hat{f}$ is \emph{optimal in the minimax sense} over the class $\mathcal{F}$ if there exists a positive sequence $(\rho_n)_n$ and positive constants $c$ and $C$ such that 
$$ \sup_{f\in \mathcal{F}}\esp{\ell(\hat{f},f)}\leq C\rho_n, 
\qquad\text{and}\qquad
\inf_{\tilde{f}} \sup_{f\in \mathcal{F}}\esp{\ell(\tilde{f},f)}\geq c\rho_n,$$
where the infimum is taken over all estimators $\tilde{f}$ of $f$ based on $(X_1, \dots, X_n)$. 
Such sequence $(\rho_n)_{n\geq 1}$ is called the minimax rate over the class $\mathcal{F}$. 
Moreover, if the estimator $\hat{f}$ does not depend on the unknown regularity of $f$, it is said to be \emph{adaptive in the minimax sense} over the class $\mathcal{F}$.
The most commonly used loss functions are defined from $\L_p$ norms. 
The literature on nonparametric density estimation is abundant.  
We can refer to the book of \cite{tsybakov2009introduction} for a review of optimal results in different contexts. 
This topic remains widely studied in more intricate contexts such as deconvolution, non independent observations, and, more recently, data privacy constraints as considered in this paper.


\paragraph{Differential privacy.}
The collection of a large amount of personal data requires new tools to protect sensitive information concerning individuals. 
To guarantee the privacy of each individual, the development of 
mechanisms to be applied to data bases has become crucial. 
Differential privacy has been widely adopted, since the seminal papers by \cite{dwork2006calibrating, dwork2006our}, as it provides rigorous privacy guarantees. The original definition corresponds to the notion of \emph{global differential privacy}, where the entire original data set $(X_1, \ldots, X_n)$ is privatized into an output that preserves the privacy of the $n$ individuals and is used for further statistical analyses. This means that the $n$ data holders share confidence with a common curator who has access to the whole sample $(X_1, \ldots, X_n)$.  In this paper, we consider a stronger notion of privacy, that is called \emph{local differential privacy}, where the individuals successively generate a private view $Z_i$ of their original data $X_i$, possibly taking into account the previous privatized data $Z_1,\ldots,Z_{i-1}$.
More precisely, $Z_i$ is a stochastic transformation of $X_i$ by the channel $Q_i$: 
the $i$-th individual generates $Z_i \sim Q_i(. | X_i=x_i,Z_1=z_1, \ldots, Z_{i-1}=z_{i-1})$. 
For a given positive number $ \alpha$, we say that the sequence of channels $(Q_i)_{i=1, \ldots, n}$ provides $\alpha$-local differentially private ($\alpha$-LDP) views of $(X_1, \ldots, X_n)$ if for all $1\leq i\leq n$, $z_1,\ldots,z_{i-1}$ in $\mathcal{Z}$, and $x_i, x_i'$ in $\mathcal X$, 
\begin{equation}
	\label{eq:localPrivacy}
	\sup_{B \in \mathcal B,   } \frac{Q_i( B| X_i=x_i,Z_1=z_1, \ldots, Z_{i-1}=z_{i-1})}{Q_i(B|X_i=x_i',Z_1=z_1, \ldots, Z_{i-1}=z_{i-1})} \leq e^\alpha,
\end{equation}
where $(\mathcal X , \mathcal A)$ and $(\mathcal Z , \mathcal B)$ respectively denote the measure spaces of the original data and of the privatized data.
Let us mention the particular case of \emph{non-interactive} local differential privacy where for all $1\leq i \leq n$, given that $X_i=x_i$, $Z_i$ is generated independently of the others, namely $Z_i \sim Q_i(. | X_i=x_i)$. In this case, the sequence of channels $(Q_i)_{i=1, \ldots, n}$ provides non-interactive $\alpha$-LDP views of $(X_1, \ldots, X_n)$ if for all $1\leq i\leq n$, and $x_i, x_i'$ in $\mathcal X$, 
\begin{equation}
	\label{eq:localPrivacyNonInter}
	\sup_{B \in \mathcal B,   } \frac{Q_i( B| X_i=x_i)}{Q_i(B|X_i=x_i')} \leq e^\alpha.
\end{equation}

In statistical inference under privacy constraints, we do not have access to the original data. Therefore, the statistical performances are necessarily deteriorated. In order to quantify this gap, the literature is prolific since the early work by  \cite{wasserman2010statistical,smith2011privacy}. The first minimax rates for estimation problems under differential privacy conditions are established by \cite{duchi2013local, duchi2013localnips, duchi2018minimax} and more recently by \cite{cai2021cost}. 
For instance, \cite{butucea2020local} propose adaptive estimators of the density over Besov ellipsoids with respect to the $\L_p$ norm, based on wavelet thresholding. \cite{sart2023density} constructs adaptive density estimators under local differential privacy with respect to the Hellinger losses. \cite{Rohde2020} deal with the estimation of functionals. \cite{butucea2023interactive} study the effect of interactive versus non interactive privacy  mechanisms on the estimation of quadratic functionals, while 
\cite{Lam2022, Dubois2023} consider goodness-of-fit testing problems.

In this paper, we are interested in adaptive minimax density estimation under local differential privacy, with respect to adversarial losses. Let us now introduce the adversarial framework.


\paragraph{Adversarial framework.} 

Beyond density estimation, generative models have seen considerable growth in recent years. Generative models aim to reproduce the sampling behaviour of a target distribution, rather than explicitly fitting a density function. 
Specifically, machine learning has made significant empirical progress in generative modeling, using such tools as generative adversarial networks (GANs) introduced by \cite{goodfellow2014GAN}. GANs provide a flexible framework for sampling from unknown distributions, and have become a standard tool among practitioners. From a practical point of view, numerous improvements have been made to GANs, enabling them to achieve state-of-the-art performance in various data generation tasks. 
We refer to the work of \cite{gui2020review} for a review and all references therein. 
The empirical success of GANs motivated many researchers to analyze their theoretical properties; among others, we mention the work of \cite{biau2020theoretical,liang2021well,schreuder2021statistical,asatryan2023convenient,puchkin2024rates}.

The rise of GANs has brought with it a focus on losses that can drive generative models, such as integral probability metrics (IPMs) as considered by \cite{muller1997integral,dziugaite2015training,liu2017approximation,bottou2018geometrical}. IPMs are types of distance between probability distributions, defined by the ability of a class of functions to distinguish between two distributions. More precisely, given a class $\mathcal{G}$ of real-valued functions defined on $\mathcal X $,
we define the distance between two probability density functions $f_1$ and $f_2$ by
\begin{equation}\label{eq:adversloss} 
	d_{\mathcal{G}}(f_1,f_2) = \sup_{g\in\mathcal{G}}\int_{\mathcal X } (f_1-f_2) g. 
\end{equation}
By choosing different sets $\mathcal{G}$, this framework can express a multitude of commonly used measures. We refer to it as the \emph{adversarial framework}. The name adversarial comes from the idea that $\mathcal{G}$ can be viewed as a \emph{discriminator} able to distinguish between the two probability measures carried by $f_1$ and $f_2$.
For example, choosing the 1-Lipschitz functions equipped with the Wasserstein metric leads to the Wasserstein GANs used by \cite{arjovsky2017wasserstein}, whereas choosing $\mathcal{G}$ as a Sobolev space leads to Sobolev GANs as considered by \cite{mroueh2018sobolev}.

In a pre-published version from 2017, \cite{liang2021well} is the first to formalize nonparametric estimation under the adversarial framework, and to prove upper bounds for Sobolev GANs. Following this work, \cite{singh2018nonparametric} succeed in improving this upper bound. Finally, \cite{liang2021well} obtain an optimal minimax rate for Sobolev GANs. 
The study of the adversarial framework is today a very active field of research, with variations on the choice of regularity classes, or by generalizing the choice of the metric. To name a few, \cite{singh2018minimax,bai2019approximability,weed2019estimation,chen2022distribution} focus on Wasserstein distance;
\cite{uppal2019nonparametric} on Besov spaces; \cite{luise2020generalization} on optimal transport-based loss functions. In this work, we focus on private density estimation using adversarial losses as defined in Equation \eqref{eq:adversloss}, where the discriminator $\mathcal{G}$ is a Sobolev ball. This can be seen as a first step towards understanding the performances of Sobolev GANs under privacy.


\paragraph{Minimax rate in the adversarial framework under privacy constraints.}
In this paper, our aim is to estimate the density of $X_1,\ldots,X_n$ under privacy constraints. 
Hence, the definition of the minimax rate needs to take into account the choice of the set of all privacy mechanisms that are considered.
More precisely, let $\mathcal{Q}_\alpha$ denote the set of all $\alpha$-LDP mechanisms, that is composed of all mechanisms $Q =(Q_1, \ldots, Q_n)$ satisfying Equation \eqref{eq:localPrivacy}. The minimax rate over the regularity set $\mathcal{F}$ with respect to the adversarial loss $d_\mathcal{G}$ under $\alpha$-local differential privacy is defined by
\begin{equation}\label{eq:defrhon}
	\rho_n(\mathcal{F},\mathcal{G},\mathcal{Q}_\alpha)
	=\inf_{Q\in\mathcal{Q}_\alpha}\inf_{\tilde{f}}\sup_{f\in\mathcal{F}} \esp{d_\mathcal{G}\pa{\tilde{f},f}},
\end{equation}
where the infimum over $\tilde{f}$ is taken among all estimators of $f$ based on privatized data $(Z_1,\ldots,Z_n)$ obtained from $Q$. 
Consider an $\alpha$-LDP mechanism $Q$ and an estimator $\hat{f}$ based on data privatized from $Q$. 
We say that the couple $(Q,\hat{f})$ is \emph{minimax optimal} if there exists a constant $C$ (w.r.t. $n$) such that 
$$
\sup_{f\in\mathcal{F}} \esp{d_{\mathcal{G}}(\hat{f},f)} \leq C \rho_n(\mathcal{F},\mathcal{G},\mathcal{Q}_\alpha).
$$
In this paper, we consider Sobolev spaces as regularity sets, either $\mathcal{F}$ for the density or $\mathcal{G}$ for the discriminator in the adversarial loss. 
Note that interactive mechanisms may lead to faster minimax rates in density estimation than non-interactive ones. However, as in the work of \cite{butucea2020local} for instance, we show that it is not the case in our framework. 
Indeed, one of the main contributions of this work is the introduction of a new non-interactive mechanism, which we refer to as the  \emph{Coordinate block privacy mechanism}, which achieves the minimax optimal rates with projection estimators on the Fourier basis of $\L_2([0,1]^d)$. 

\paragraph{}
The structure of this paper is as follows. 
Section \ref{sect:results} consists of all the results in the isotropic case. The definition of the Sobolev balls is detailed in Section \ref{sect:Fourier}. Section \ref{sect:lower_bound} is devoted to the lower bound of the minimax rate. Section \ref{sect:upper_bound} presents the Coordinate block privacy mechanism together with our private projection estimator which are minimax optimal. 
In Section \ref{sect:adaptupper_bound}, we construct an adaptive procedure based on the method of \cite{MR2850214} which is minimax optimal up to logarithmic terms.
In Section \ref{sect:anisotrop}, we present results in an anisotropic framework. 
In Section \ref{sect:proofs}, we prove the main results. The proofs are detailed in the isotropic case to avoid making the notations too cumbersome. Complementary proofs can be found in Appendix \ref{sect:complproofs}. Finally, the necessity of the Coordinate block privacy mechanism to achieve minimax optimal rates is explained in Appendix \ref{app:globalnotminimax}. \\

In the following, $c$, $C$ and $c_{a,b,\ldots}$, $C_{a,b,\ldots}$ denote respectively universal constants and constants depending only on $a,b,\ldots$, that may vary from line to line. Moreover, the notation $u_n \asymp v_n$ indicates that there exist constants $c$ and $C$ such that for all $n$, $c\leq u_n/v_n \leq C.$

\section{Results over isotropic Sobolev balls}\label{sect:results}

Before stating our main results, let us define the regularity spaces considered in this section, namely isotropic Sobolev balls.

\subsection{Fourier basis and Sobolev spaces}
\label{sect:Fourier}

Consider the $d$-dimensional Fourier basis of $\L_2([0,1]^d)$ obtained by taking the tensor product of the one-dimensional Fourier basis. 
More precisely, recall the one-dimensional Fourier basis defined for all $t$ in $[0,1]$ by
$$ 
\varphi_1(t) = 1, \qquad \text{and}\qquad\forall j\in\N^*, \quad
\left\{\begin{array}{lcl}
	\varphi_{2j}(t) &=& \sqrt{2} \cos(2\pi jt) \\
	\varphi_{2j+1}(t) &=& \sqrt{2} \sin(2\pi jt).
\end{array}\right.
$$
The $d$-dimensional Fourier basis is defined for all $\mi{j}=(j_1,\ldots,j_d)$ in $(\N^*)^d$ and for all $x = (x_1,\ldots,x_d)$ in $[0,1]^d$ by 
$$
\varphi_{\mi{j}}(x) = \prod_{m=1}^{d} \varphi_{j_m}(x_m).
$$

Then, the family $\{\varphi_{\mi{j}}\}_{\mi{j}\in(\N^*)^d}$ is an orthonormal basis of $\L_2([0,1]^d)$. In particular, any function $f$ in $\L_2([0,1]^d)$ can be uniquely decomposed as 
$$
f=\sum_{\mi{j}\in(\N^*)^d} \theta_{\mi{j}}(f)\varphi_{\mi{j}}, 
\qquad \text{where}\qquad 
\theta_{\mi{j}}(f) = \int_{[0,1]^d} f(x) \varphi_{\mi{j}}(x) dx.
$$

A natural choice for regularity spaces are Sobolev balls since they can be simply characterized in terms of the Fourier decomposition. The isotropic Sobolev balls are defined for all smoothness parameter $\beta>0$ and radius $R>0$ by 
\begin{equation}\label{eq:defSobolev}
	\mathcal{W}^{\beta}(R) = \ac{f\in\L_2([0,1]^d) \; ; \; 
		\sum_{\mi{j}\in(\N^*)^d} \pa{j_1^{2\beta} + \ldots + j_d^{2\beta}}\theta_{\mi{j}}^2(f)\leq R^2}.
\end{equation}

Note that if $f$ is a density, then 
$$
\sum_{\mi{j}\in(\N^*)^d} \pa{j_1^{2\beta} + \ldots + j_d^{2\beta}}\theta_{\mi{j}}^2(f) 
\geq d\ \theta_{(1,1,\ldots,1)}^2(f) = d.$$
Hence the Sobolev ball $\mathcal{W}^{\beta}(R)$ may contain densities only if $R^2 \geq d$. 
In this section, both densities and discriminant functions will belong to such isotropic Sobolev balls.

\subsection{Lower bound over Sobolev balls}
\label{sect:lower_bound}

In this section, we bound from below the minimax rate over Sobolev balls in the adversarial framework under privacy constraints. 

\begin{theorem}
	\label{thm:lowerbound}
	Consider a sample size $n$, a privacy level $\alpha$ in $(0,A]$ such that $n\alpha^2\geq 1$, integer smoothness parameters $\beta$ and $\delta$ in $\N^*$, and a positive radius $R$ such that $R^2>d$. 
	Then, the minimax rate defined in \eqref{eq:defrhon} is bounded from below as follows:
	$$
	\rho_n(\mathcal{W}^{\beta}(R),\mathcal{W}^{\delta}(1),\mathcal{Q}_\alpha) \geq
	c\max\ac{(n\alpha^2)^{-\frac{\beta+\delta}{2\beta+2d}},(n\alpha^2)^{-1/2}}, 
	$$
	where $c=c_{A,d,\delta,\beta,R} $ denotes a positive constant depending on $A$, $d$, $\delta$, $\beta$, and $R$. 
\end{theorem}

The proof of Theorem \ref{thm:lowerbound} is detailed in Section \ref{sect:proofthmlowerbound}. Depending on the value of the smoothness parameter $\delta$ for the adversarial distance, the rate is not the same. In particular, if $\delta>d$, we obtain a near-parametric rate, while for $\delta<d$, we see the influence of the regularity of the discriminator class $\mathcal{G}=\mathcal{W}^{\delta}(1)$. Note that these rates are coherent with the ones achieved by \cite{liang2021well} over slightly different Sobolev spaces, with respect to adversarial losses in dimension $d$, without privacy constraints, that are $\max\{n^{-\frac{\beta+\delta}{2\beta+d}}, n^{-1/2}\}$. One may notice the effect of privatization, that replaces $n$ by $n\alpha^2$ and transforms the dimension $d$ into $2d$. 
These effects can also be observed by \cite{duchi2013localnips}, over Sobolev spaces in dimension 1, with respect to the $\L_2$ loss, where the rate $n^{-2\beta/(2\beta+1)}$ without confidentiality becomes $(n\alpha^2)^{-2\beta/(2\beta+2)}$ under privacy constraints. 

The lower bound is proved only for integer regularity parameters $\beta$ and $\delta$ which is standard in the literature. Indeed, the proof is based on the construction of parametric spaces, namely $\mathcal{F}^{\beta}$ and $\mathcal{D}^\delta$, which relies on the characterization of Sobolev spaces based on the $\L_2$ norm of the partial derivatives (see Lemma \ref{lm:derivtocoefSobolev} in Appendix \ref{sect:sobolevcharact}).

\subsection{Optimal local differential private projection estimator}
\label{sect:upper_bound}

In the usual \emph{non private} setting, projection estimators have been widely studied. Over isotropic Sobolev balls $\mathcal{W}^\beta(R)$, as considered in this section, where the smoothness parameter $\beta$ is the same for each dimension, the density $f$ of the original data $(X_1,\ldots,X_n)$ is well-approximated by its orthogonal projection $f_J$ on $\operatorname{Span}\pa{\varphi_{\mi{j}}, \mi{j}\in\{1,\ldots,J\}^d }$, namely
\begin{equation}
	\label{eq:projorthof}
	f_{J}=\sum_{\mi{j}\in\{1,\ldots,J\}^d} \theta_{\mi{j}}(f)\varphi_{\mi{j}}, 
\end{equation}
for a well-chosen $J$. 
Hence, if we had access to the original data,  we could then estimate the density $f$ by the empirical counterpart of $f_J$, that is
$$
\tilde{f}_J = \sum_{\mi{j}\in\ac{1,2, \ldots,J}^d} \tilde{\theta}_{\mi{j}} \varphi_{\mi{j}}, 
\qquad \text{where}\qquad \forall \mi{j}\in\ac{1,2, \ldots,J}^d,\ 
\tilde\theta_{\mi{j}} = \frac{1}{n}\sum_{i=1}^n \varphi_{\mi{j}}(X_i),
$$
and where the integer $J$ needs to be calibrated. 
Yet, under privacy constraints, we are not allowed to use the original dataset $(X_i)_{1\leq i\leq n}$, only a privatized version of the data. In Section \ref{sect:priv_machanism}, we first define a general mechanism $Q$ that provides non-interactive $\alpha$-local differential private data, by independently disturbing blocks of coordinates $(\varphi_{\mj}(X_i))_{\mj}$, independently for each individual. In Section \ref{sect:projEstim}, we apply this new mechanism with dyadic blocks of $\{1,\ldots,J\}^d$, and denote it $Q^{(J)}$. We then introduce the private projection estimator $\hat{f}_J$. Finally, we prove that the couple $(Q^{(J)},\hat{f}_J)$ achieves optimal rates.

\subsubsection{Coordinate block privacy mechanism}
\label{sect:priv_machanism}

This mechanism is strongly inspired by the one proposed by \cite{duchi2018minimax} and adapted by \cite{butucea2023phase}. 
The main difference is that, instead of privatizing all the coefficients at once, as done in the \emph{Coordinate global privacy mechanism} by \cite{butucea2023phase}, we privatize independently blocks of coordinates. The other difference is that we adapt the mechanism such that it can be applied to blocks of size 2. \\

For an individual $X$, we want to generate a private version $Z$ of $X$ by randomly disturbing $\varphi(X) = (\varphi_{\mj}(X))_{\mj\in\mcJ},$ where $\mcJ$ is a finite set of indices. Assume that we are given blocks $\{\mathcal{J}_{\mi{\ell}}, \mi{\ell}\in\mathcal{L}\}$ that form a partition of $\mathcal{J}$, and denote $d_{\ml}$ the cardinal of $\mcJ_{\ml}$. Fix a global privacy level $\alpha>0$, and associate to each block a privacy level $\alpha_{\ml}$ such that $\sum_{\mi{\ell}\in\mcL}\alpha_{\mi{\ell}} = \alpha$. 
Then, we define the mechanism $Q_1$ that generates a private version $Z$ of an individual $X$ using the steps below. 
\begin{enumerate}
	\item Compute the coordinates of $\varphi(X)$ in each block, denoted
	$$
	\varphi(X) = \pa{\varphi_{[\ml]}(X)}_{\ml\in\mcL} 
	\quad \text{with} \quad 
	\varphi_{[\ml]}(X) = \pa{\varphi_{\mj}(X)}_{\mj\in\mcJ_{\ml}} \in \cro{-B_0,B_0}^{d_{\mj}}, 
	$$ 
	where $B_0=\max_{\mi{j}\in\mcJ}\norm{\varphi_{\mi{j}}}_{\infty}.$ Note that, for the Fourier basis as considered in this paper, $B_0=2^{d/2}$. If this mechanism is applied to more general $\varphi_{\mj}$ that are not bounded, then one can apply the censoring operator of \cite{butucea2023phase}. 
	\item Independently for all $\mi{\ell}$ in $\mcL$, apply the steps below. 
	\begin{enumerate}
		\item Draw at random vectors $\tilde{V}_{[\mi{\ell}]}$ in $\{-B_0,B_0\}^{d_{\mi{\ell}}}$ with independent coordinates such that
		$$
		\forall\mi{j}\in \mathcal{J}_{\mi{\ell}},\quad \tilde{V}_{\mi{j}} = \left\{\begin{array}{ll}
			B_0 & \text{with probability }\frac{1}{2} + \frac{\varphi_{\mi{j}}(X)}{2B_0}, \\ 
			-B_0 & \text{otherwise}.
		\end{array}\right.
		$$ 
		\item Compute 
		\begin{equation}\label{eq:defBdalpha}
			B_{d_{\ml}}\pa{\alpha_{\ml}} = B_0 \frac{e^{\alpha_{\ml}}+1}{e^{\alpha_{\ml}}-1} \Gamma_{d_{\ml}}, 
			\qquad \text{with}\qquad 
			\frac{1}{\Gamma_{d_{\ml}}} =  \frac{1}{2^{d_{\ml}-1}}\binom{d_{\ml}-1}{\lfloor\frac{d_{\ml}-1}{2}\rfloor}. 
		\end{equation}
		\item Conditionally on $\tilde{V}_{[\ml]}=\tilde{v}_{[\ml]}$, and independently on $X$, 
		draw at random 
		$$
		Y_{\mi{\ell}} \sim \mathcal{B}(p_{\ml}),
		\quad \text{where}\quad 
		p_{\ml} = \left\{\begin{array}{ll} 
			\displaystyle 1 & \text{if $d_{\ml}$ is odd,} \\ 
			\displaystyle 1 - \frac{1}{2^{d_{\ml}}} \binom{d_{\ml}}{d_{\ml}/2} & \text{if $d_{\ml}$ is even.} 
		\end{array}\right.
		$$
		\begin{enumerate}
			\item 
			If $Y_{\ml}=1$, then generate 
			$$
			T_{\mi{\ell}}\sim\mathcal{B}(\pi_{\alpha_{\mi{\ell}}}),
			\quad \text{where}\quad \pi_{\alpha_{\mi{\ell}}} = \frac{e^{\alpha_{\mi{\ell}}}}{1+e^{\alpha_{\mi{\ell}}}},
			$$ 
			and
			$$
			Z_{[\mi{\ell}]} \sim \left\{\begin{array}{ll}
				\mathcal{U}\pa{\ac{z\in\{\pm B_{d_{\mi{\ell}}}(\alpha_{\mi{\ell}})\}^{d_{\mi{\ell}}} \; ; \; \langle z,\tilde{v}_{[\ml]} \rangle >0}} 
				& \text{if }T_{\mi{\ell}}=1,\\
				\mathcal{U}\pa{\ac{z\in\{\pm B_{d_{\mi{\ell}}}(\alpha_{\mi{\ell}})\}^{d_{\mi{\ell}}} \; ; \; \langle z,\tilde{v}_{[\ml]} \rangle <0}} 
				& \text{if }T_{\mi{\ell}}=0.
			\end{array}\right.
			$$
			\item If $Y_{\ml}=0$, then generate 
			$$
			Z_{[\mi{\ell}]} \sim \mathcal{U}\pa{\ac{z\in\{\pm B_{d_{\mi{\ell}}}(\alpha_{\mi{\ell}})\}^{d_{\mi{\ell}}} \; ; \; \langle z,\tilde{v}_{[\ml]} \rangle = 0}}.
			$$
		\end{enumerate}
	\end{enumerate}
\end{enumerate}

First note that if $d_{\ml}$ is odd, this mechanism applied to a single block coincides with both the mechanisms of \cite{duchi2018minimax} and \cite{butucea2023phase}. Indeed, in this case, $p_{\ml}=1$ and thus $Y_{\ml}=0$ never occurs. Besides, for all $z$ in $\{\pm B_{d_{\mi{\ell}}}(\alpha_{\mi{\ell}})\}^{d_{\mi{\ell}}}$ and $\tilde{v}$ in $\{\pm B_0\}^{d_{\mi{\ell}}}$, $\langle z,\tilde{v} \rangle = 0$ is impossible. 
Indeed, the scalar product $\langle z,\tilde{v} \rangle$ is proportional to the difference between the number of coordinates in $z$ and $\tilde{v}$ of same sign and the number of ones of opposite sign. It can thus be equal to zero only if the number of coordinates is even. 

Besides from the use of blocks, the novelty comes from the fact that we deal separately the cases where the scalar product $\langle z,\tilde{v}_{[\ml]} \rangle$ equals zero, which happens with probability $p_{\ml}$ chosen such that
$$
p_{\ml} 
= 1 - \frac{\abs{\ac{z\in\{\pm B_{d_{\mi{\ell}}}(\alpha_{\mi{\ell}})\}^{d_{\mi{\ell}}} \; ; \; \langle z,\tilde{v}_{[\ml]} \rangle =0}}}{2^{d_{\ml}}}, 
$$
as proved in Equation \eqref{eq:calculpell}. 

\begin{proposition}
	\label{prop:CBmechanismalphaLDP}
	Consider the mechanism $Q_1$ as defined above, and $Q$ the mechanism which generates independently for each individual $i$, $Z_i = (Z_{i,\mj})_{\mj\in\mcJ} \sim Q_1(\cdot|X_i)$. Then, $Q$ provides non-interactive $\alpha$-LDP views of $X_1,\ldots,X_n$, as defined in \eqref{eq:localPrivacyNonInter}.
\end{proposition}

\begin{proposition}
	\label{prop:CBmechanismespvar}
	Consider the mechanism $Q_1$ as defined above and $Z = (Z_{\mj})_{\mj\in\mcJ} \sim Q_1(\cdot|X)$. 
	\begin{enumerate}
		\item Then, for all $\mj$ in $\mcJ$, 
		$$\esp{Z_{\mi{j}}\middle|X} = \varphi_{\mi{j}}(X).$$
		\item Moreover, for all $\mi{\ell}$ in $\mcL$, and $\mi{j}$ in $\mathcal{J}_{\mi{\ell}}$,  
		$$
		\abs{Z_{\mi{j}}} \leq 2 \xi_A B_0 \frac{\sqrt{d_{\mi{\ell}}}}{\alpha_{\mi{\ell}}}
		\quad \text{thus}\quad 
		\var{Z_{\mi{j}}} \leq (2 \xi_A B_0)^2\frac{d_{\mi{\ell}}}{\alpha_{\mi{\ell}}^2}, 
		$$
		where $\xi_A = A(e^A+1)/(e^A-1)$. 
	\end{enumerate}
\end{proposition}

The proofs are quite similar to the ones of \cite{butucea2023phase} and may be found in Sections \ref{sect:proof:prop:CBmechanismalphaLDP} and \ref{sect:proof:prop:CBmechanismespvar}. 
A first natural idea was to consider the \emph{Coordinate global privacy mechanism} of \cite{duchi2018minimax} and \cite{butucea2023phase}, with $(\varphi_{\mj})_{\mj}$ the Fourier basis, and $\mcJ=\{1,\ldots,J\}^d$. Yet, one may prove that the variance of each $Z_{i,\mi{j}}$ is then of order $J^d/\alpha^2$ and this allows to achieve optimal rates only for $\delta<d/2$. This  point is explained in detail in Appendix \ref{app:globalnotminimax}.
With this new mechanism, the variance of $Z_{i,\mi{j}}$ is controlled by the size $d_{\mi{\ell}}$ of the corresponding block, and it turns out that it makes it possible to achieve the optimal minimax rate. 
Finally, note that, in view of developing an adaptive estimator, the constant in the upper bound of the variance needs to be computable, which is the case here. 
Moreover, considering $\xi_A$ provides a much finer upper bound than the more usual $e^A+1$, as it is equivalent to $A$ when $A\to+\infty$.

\subsubsection{Private projection estimator}
\label{sect:projEstim}

Let us now introduce dyadic blocks. Consider $J=2^{L+1}-1$, where $L\geq 0$, and define for all $\mi{\ell}=(\ell_1,\ldots,\ell_d)$ in $\N^d$, 
$$
\mathcal{J}_{\mi{\ell}} = \prod_{m=1}^d \ac{2^{\ell_m},\ldots,2^{\ell_m+1}-1},
$$ 
such that the blocks $\mathcal{J}_{\mi{\ell}}$ with $\mi{\ell}$ in $\ac{0,\ldots,L}^d$ form a partition of $\mcJ=\ac{1,2, \ldots,J}^d.$
Denote $d_{\mi{\ell}}$ the cardinal of $\mathcal{J}_{\mi{\ell}}$, that is 
\begin{equation}\label{eq:defdell}
	d_{\mi{\ell}} = |\mathcal{J}_{\mi{\ell}}| = \prod_{m=1}^d 2^{\ell_m}.
\end{equation}
To each block $\mathcal{J}_{\mi{\ell}}$, associate a privacy level $\alpha_{\mi{\ell}}$ that will be calibrated later (see Equation  \eqref{eq:defalphaell}), and such that $\sum_{\mi{\ell}\in\ac{0,1,\ldots,L}^d}\alpha_{\mi{\ell}} = \alpha$. 
Note that, when dealing with dyadic blocks, as done in this paper, the only $\mi{\ell}$ for which $d_{\mi{\ell}}$ is odd is the case $\mi{\ell}=(0,0,\ldots,0)$. 
Let us define the \emph{Coordinate dyadic block privacy mechanism} $Q^{(J)}$ as the non-interactive mechanism described in Section \ref{sect:priv_machanism} relying on the dyadic blocks of $\{1,\ldots,J\}^d$ defined above. 
\\

Given the observation of the private views $Z_1,\ldots,Z_n$ of $X_1,\ldots,X_n$ generated using $Q^{(J)}$, we aim at estimating the density $f$ of the $X_i$'s. 
As done by \cite{duchi2013localnips} in dimension one, we estimate the coefficients $\theta_{\mi{j}}(f)$ in the Fourier expansion of $f$ for all $\mi{j}$ in $\ac{1,2,\ldots,J}^d$ by
$$
\hat\theta_{\mi{j}} = \frac{1}{n} \sum_{i=1}^n Z_{i,\mi{j}}.
$$
Then, we estimate the density $f$ by 
\begin{equation}
	\label{eq:privateprojestim}
	\hat{f}_J = \sum_{\mi{j}\in\{1,2,\ldots,J\}^d} \hat\theta_{\mi{j}} \varphi_{\mi{j}},
\end{equation}
for a well calibrated $J$. 
In particular, Proposition \ref{prop:CBmechanismespvar} implies that for all $\mi{j}$, $\hat{\theta}_{\mi{j}}$ is an unbiased estimator of $\theta_{\mi{j}}(f)$. 
Thus, $\hat{f}_J$ is an unbiased estimator of the orthogonal projection $f_J$ defined in \eqref{eq:projorthof}. 
We then obtain a bias-variance decomposition 
\begin{equation}\label{eq:biaisvardecomp}
	\esp{d_{\mathcal{W}^{\delta}(1)}\pa{\hat{f}_J,f}}\leq \esp{d_{\mathcal{W}^{\delta}(1)}\pa{\hat{f}_J,f_J}} + d_{\mathcal{W}^{\delta}(1)}\pa{f_J,f}. 
\end{equation}

The first term in the right hand side of \eqref{eq:biaisvardecomp} is a variance term and corresponds to the stochastic error. It is controlled using properties on the privacy mechanism, as below. 

\begin{proposition}\label{prop:controlvar}
	Consider a sample size $n$, a privacy level $\alpha$ in $(0,A]$ such that $n\alpha^2> 1$, and a smoothness parameter $\delta$ in $\R_+^*$. Fix $J = 2^{L+1}-1$ for some $L$ in $\N$. 
	Consider the Coordinate dyadic block privacy mechanism $Q^{(J)}$ and the projection estimator $\hat{f}_J$ defined in Equation \eqref{eq:privateprojestim}. Then,  
	\begin{equation}\label{eq:majvar}
		\esp{d_{\mathcal{W}^{\delta}(1)}\pa{\hat{f}_J,f_J}} \leq \tau_{A,d} \Sigma_J,
	\end{equation}
	where 
	\begin{equation}\label{eq:defsigmaell}
		\displaystyle\tau_{A,d}=2\sqrt{\frac{2^{d}}{d}}\frac{A(e^A+1)}{(e^A-1)}, \quad 
		\Sigma_J = \sum_{\mi{\ell} \in \{0,\ldots,L\}^d} \sigma_{\mi{\ell}},
		\quad \text{and}\quad 
		\sigma_{\mi{\ell}} = \frac{1}{\sqrt{n}}\frac{d_{\mi{\ell}}^{1-\delta/d}}{\alpha_{\mi{\ell}}}.
	\end{equation}
\end{proposition}
The proof is detailed in Section \ref{sect:proof:prop:controlvar}. In particular, solving the constrained-optimization problem which consists in 
$$
\min_{(\alpha_{\mi{\ell}})_{\mi{\ell}}} \ac{\sum_{\mi{\ell}\in \{0,1,\ldots,L\}^d} \frac{d_{\mi{\ell}}^{1-\delta/d}}{\alpha_{\mi{\ell}}}} 
\qquad \text{subject to}\qquad \sum_{\mi{\ell}\in \{0,1,\ldots,L\}^d} \alpha_{\mi{\ell}} = \alpha,
$$ 
by the method of Lagrangian multipliers leads to the choice of the block privacy levels 
\begin{equation}\label{eq:defalphaell}
	\alpha_{\mi{\ell}} = \frac{\alpha}{S_{d,\delta}(J)} d_{\mi{\ell}}^{(1-\delta/d)/2} 
	\qquad \text{where}\qquad 
	S_{d,\delta}(J) = \sum_{\mi{\ell}\in\{0,1,\ldots,L\}^d} d_{\mi{\ell}}^{(1-\delta/d)/2}. 
\end{equation}
Such block privacy levels $\alpha_{\mi{\ell}}$ minimize the variance term in the bias-variance decomposition \eqref{eq:biaisvardecomp}.
In particular, one may notice that they are proportional to $d_{\mi{\ell}}^{(1-\delta/d)/2}$. Hence, if the discriminator regularity $\delta<d$, then one needs to consider larger privacy levels (and thus less confidentiality) for larger blocks, whereas if $\delta>d$, smaller privacy levels (and thus more confidentiality) need to be taken for larger blocks. 
Finally, if $\delta=d$, then $\alpha_{\mi{\ell}} = \alpha/(L+1)^d,$ and all blocks have the same privacy level. Note that because of the constraint, the block privacy levels $\alpha_{\mi{\ell}}$ depend on $J$, and thus, so does the privacy mechanism $Q^{(J)}$.

For this particular choice of privacy levels, we obtain that the upper-bound of the variance term is controlled by
\begin{equation}\label{eq:defSigmaJ}
	\Sigma_J= \frac{S_{d,\delta}(J)^2}{\sqrt{n\alpha^2}}.
\end{equation}

The second term in the right hand side of \eqref{eq:biaisvardecomp} is a bias term which corresponds to the approximation error. It is controlled using the regularity of the density $f$. 
\begin{lemma}\label{lm:controlbias}
	Consider positive smoothness parameters $\beta$ and $\delta$, and a positive radius $R$. 
	If $f$ belongs to $\mathcal{W}^{\beta}(R)$, then 
	$$d_{\mathcal{W}^{\delta}(1)}\pa{f_J,f} \leq R J^{-(\beta+\delta)}.$$
\end{lemma}

The proof of Lemma \ref{lm:controlbias} is detailed in Section \ref{sect:proof:lm:controlbias}.
Finally, a trade-off between the bias and the variance terms leads to Theorem \ref{thm:upper}.
\begin{theorem}\label{thm:upper}
	Consider a sample size $n$, a privacy level $\alpha$ in $(0,A]$ such that $n\alpha^2> 1$, positive smoothness parameters $\beta$ and $\delta$, and a positive radius $R$ such that $R^2\geq d$. Let $J=2^{L+1}-1$ for some $L$ in $\N$.
	Consider the Coordinate dyadic block privacy mechanism $Q^{(J)}$, where the block privacy levels $(\alpha_{\mi{\ell}})_{\mi{\ell}}$ are defined in \eqref{eq:defalphaell}.
	Then, the projection estimator $\hat{f}_J$ defined in \eqref{eq:privateprojestim} satisfies 
	$$
	\sup_{f\in\mathcal{W}^{\beta}(R)} \esp{d_{\mathcal{W}^{\delta}(1)}\pa{\hat{f}_J,f}} \leq
	\left\{ \begin{array}{ll}
		C (n\alpha^2)^{-\frac{\beta+\delta}{2\beta+2d}} & \text{if $\delta<d$ and }J\asymp (n\alpha^2)^{\frac{1}{2\beta+2d}} \\
		C \pa{\frac{n\alpha^2}{\cro{\log(n\alpha^2)}^{4d}}}^{-1/2} & \text{if $\delta=d$ and }J\asymp\pa{\frac{n\alpha^2}{\cro{\log(n\alpha^2)}^{4d}}}^{\frac{1}{2\beta+2\delta}}\\
		C (n\alpha^2)^{-1/2} & \text{if $\delta>d$ and }J\asymp(n\alpha^2)^{\frac{1}{2\beta+2\delta}},
	\end{array}\right.
	$$
	where $C=C_{A,d,\delta,\beta,R} $ denotes a constant depending on $A$, $d$, $\delta$, $\beta$, and $R$. 
\end{theorem}

The proof is detailed in Section \ref{sect:proof:thm:upper}. Unlike the case of the lower bound, the result also holds for non-integer parameters $\beta$ and $\delta$. 
One may see with Theorem \ref{thm:lowerbound} that the private projection estimator with the Coordinate dyadic block privacy mechanism is minimax optimal for all $\delta \neq d$. 
On the contrary, in the limiting case where $\delta=d$, a logarithmic term appears. In this case, we do not know if the lower bound or the upper bound are suboptimal.

Finally, as the choice of $J$ leading to optimal minimax rates depends on the unknown regularity $\beta$ of the density, this procedure is not adaptive in the minimax sense. 
The construction of an adaptive procedure is the topic of the next section.

\subsection{Adaptive local differential private projection estimator}
\label{sect:adaptupper_bound}

In this section, we propose an adaptive procedure of our method based on private data. First, note that the set of generators $\mathcal{W}^{\delta}(1)$ is chosen by the user. Then, the regularity  $\delta$ of the generator is assumed to be known from now on. We focus on the adaptation problem with respect to the regularity parameter $\beta$ of the density $f$. 

Ideally, the parameter $J$ would be chosen by minimizing the expected risk
$$
\operatorname{Crit}(J)=\esp{d_{\mathcal{W}^{\delta}(1)}\pa{\hat{f}_J,f}}
$$
for $J$ in a finite collection of models $\mathcal{M}_n^\alpha$. 
In view of the optimal choices of $J$ in Theorem \ref{thm:upper}, we consider only values of $J$ less than $n\alpha^2$. More precisely, we define the model collection by 
\begin{equation}\label{eq:defMn}
	\mathcal{M}_n^\alpha=
	\ac{J=2^{L+1}-1 \ ;\ 0 \leq L \leq \lfloor \log_2(1+n\alpha^2)\rfloor -1}. 
\end{equation}

However, this criterion is not computable in practice, since it depends on the unknown density $f$. 
We follow the method of \cite{MR2850214}, which has been introduced to select the bandwidth of a kernel estimator. Our aim is to replace the bias and variance terms in the decomposition \eqref{eq:biaisvardecomp} of the expected risk, respectively by an empirical counterpart $\hat{A}(J)$ and an upper bound $V(J)$, without making any assumption on the regularity of $f$. Then, we select $J$ that minimizes the sum of the two terms. 
This leads to 
\begin{equation}\label{eq:choixJhat}
	\hat{J}=\argmin_{J\in \mathcal{M}_n^\alpha} \widehat{\operatorname{Crit}}(J), 
	\quad \text{where} \quad 
	\widehat{\operatorname{Crit}}(J) = \hat{A}(J) + 2V(J),
\end{equation}
where the construction of $\hat{A}(J)$ and $V(J)$ is detailed below. \\

On the one hand, we define an upper bound $V(J)$ of the variance term $\Sigma_J$ introduced in Equation \eqref{eq:defsigmaell} by
\begin{equation}\label{eq:defVJ}
	V(J) = \sqrt{2}\tau_{A,d}\Sigma_J\sqrt{d\log(J) + \frac{3}{2}\log(n\alpha^2) + \log(\tau_{A,d}\Sigma_J)}, 
\end{equation}
where 
$\tau_{A,d}=2\sqrt{(2^{d}/d)}A(e^A+1)/(e^A-1)$ is also defined in \eqref{eq:defsigmaell}. Note that for any $A\geq 1$, $\tau_{A,d}\geq 2\sqrt{2}$ as $2^d/d\geq 2$. Hence, for such $A$, $V(J)$ is well-defined. Indeed, by Equation \eqref{eq:defSigmaJ}, for block privacy levels $\alpha_{\ml}$ defined in \eqref{eq:defalphaell},  
$$
d\log(J) + \frac{3}{2}\log(n\alpha^2) + \log(\tau_{A,d}\Sigma_J) 
= d\log(J) + \log(\tau_{A,d}) + \log(n\alpha^2) + 2\log\pa{S_{d,\delta}(J)} > 0, 
$$ 
as $S_{d,\delta}(J) \geq d_{(0,\ldots,0)}^{(1-\delta/d)/2} = 1$.
This choice of $V(J)$ comes from the concentration inequality below. 
\begin{lemma}\label{lm:ineq}
	Consider a sample size $n$, a privacy level $\alpha$ in $(0,A]$ such that $n\alpha^2> 1$, and $A\geq 1$. Let $\delta$ be a smoothness parameter in $\R_+^*$. Fix $J$ in $\mathcal{M}_n^\alpha$. 
	Consider the Coordinate dyadic block privacy mechanism $Q^{(J)}$ with block privacy levels $\alpha_{\ml}$ defined in \eqref{eq:defalphaell}, and the projection estimator $\hat{f}_J$ defined in  \eqref{eq:privateprojestim}. 
	Then, for all positive $t$,
	$$
	\proba{d_{\mathcal{W}^{\delta}(1)}\pa{\hat{f}_{J},f_J} \geq V(J) + t} 
	\leq \frac{1}{(n\alpha^2)^{3/2}}\times\frac{2}{\tau_{A,d}\Sigma_J}\exp\pa{\frac{-t^2}{2(\tau_{A,d}\Sigma_J)^2}}. 
	$$
\end{lemma}
The proof of Lemma \ref{lm:ineq} is detailed in Section \ref{sect:proof:lm:ineq}. \\

On the other hand, we define the empirical counterpart of the bias by
\begin{equation}
	\label{eq:defAhat}
	\hat{A} (J)= \max_{ J'\in \mathcal{M}_n^\alpha}\left\{  \left( d_{\mathcal{W}^{\delta}(1)}  
	\pa{\hat f_{J'},\hat f_{J'\wedge J}}-2V(J')\right)_{+}
	\right\},
\end{equation}
where 
$J\wedge J'=\min(J,J')$, and $(x)_{+}=\max(0,x)$ denotes the positive part of $x$. 
In Equation \eqref{eq:controlespAhatordrebiais} in the proof of Theorem \ref{thm:oraclebound} (see Section \ref{sect:proof:thm:oraclebound}), we prove that 
\begin{equation*}
	\esp{\hat{A}(J)}\leq 2\max_{\underset{J'\geq J}{J'\in \mathcal{M}_n^\alpha}}d_{\mathcal{W}^{\delta}(1)}  \pa{f, f_{J'}}+\frac{C}{n\alpha^2}.
\end{equation*}
In particular, this ensures that $\hat{A}(J)$ is of the same order as the bias. \\

Finally, the estimator $\hat{f}_{\hat{J}}$, with $\hat{J}$ is selected by minimizing the criterion defined in \eqref{eq:choixJhat}, satisfies the oracle inequality below.

\begin{theorem}[Oracle inequality]\label{thm:oraclebound}
	Consider a sample size $n$, a privacy level $\alpha$ in $(0,A]$ such that $n\alpha^2 \geq 2$ and $A\geq 1$. 
	Let $\delta$ be a positive smoothness parameter in $\R_+^*$. 
	Consider the model collection $\mathcal{M}_n^\alpha$ defined in \eqref{eq:defMn}. 
	For all $J$ in $\mathcal{M}_n^\alpha$, consider the Coordinate dyadic block privacy mechanism $Q^{(J)}$, where the block privacy levels are defined in \eqref{eq:defalphaell}, and the corresponding estimator $\hat{f}_J$ in \eqref{eq:privateprojestim}. 
	Finally, select $\hat{J}$ that minimizes the criterion $\widehat{\operatorname{Crit}}(J)$ defined by \eqref{eq:choixJhat}, \eqref{eq:defVJ} and \eqref{eq:defAhat}. 
	Then, the selected estimator $\hat{f}_{\hat{J}}$ satisfies 
	\begin{multline*}
		\esp{d_{\mathcal{W}^{\delta}(1)}\!\pa{\hat{f}_{ \hat{J}},f}} \\
		\leq 
		C\inf_{J_0\in \mathcal{M}_n^\alpha}\! \left( \esp{d_{\mathcal{W}^{\delta}(1)}\!\pa{\hat{f}_{J_0},f}} + 
		\max_{\underset{J\geq J_0}{J\in \mathcal{M}_n^\alpha}}d_{\mathcal{W}^{\delta}(1)}\!\pa{f_J, f} + V(J_0)
		\right)+\frac{C}{n\alpha^2}. 
	\end{multline*}
\end{theorem}

The proof of Theorem \ref{thm:oraclebound} is detailed in Section \ref{sect:proof:thm:oraclebound}. 
This result proves that the selected estimator $\hat{f}_{\hat{J}}$ does as well, in terms of expected risk, as the best estimator in the collection of estimators $\hat{f}_{J_0}$ for $J_0$ in $\mathcal{M}_n^\alpha$, to within terms of the same order and a negligible term. Indeed,
the term $\max_{J\in \mathcal{M}_n^\alpha, J\geq J_0}d_{\mathcal{W}^{\delta}(1)}  \pa{f_J, f}$ has the same order as the bias $d_{\mathcal{W}^{\delta}(1)}  \pa{f_{J_0}, f}$. Moreover, the term $\esp{d_{\mathcal{W}^{\delta}(1)}\pa{\hat{f}_{J_0},f}}$ is studied in Theorem \ref{thm:upper} and is smaller than the sum of the bias and variance terms. More precisely, for every $J_0\geq 1$, $\tau_{A,d} \Sigma_{J_0} \leq C_{A,d,\delta} V(J_0)$, hence
\begin{eqnarray*}
	\esp{d_{\mathcal{W}^{\delta}(1)}\pa{\hat{f}_{J_0},f}}
	\leq d_{\mathcal{W}^{\delta}(1)}  \pa{f_{J_0},f}+V(J_0).
\end{eqnarray*}
Thus, we deduce that the bias-variance decomposition of the expected risk for the selected estimator is 
\begin{equation} \label{eq:biaisvaradapt}
	\esp{d_{\mathcal{W}^{\delta}(1)}\pa{\hat{f}_{ \hat{J}},f}}
	\leq C_{A,d,\delta}\inf_{J_0\in \mathcal{M}_n^\alpha} \pa{\max_{\underset{J\geq J_0}{J\in \mathcal{M}_n^\alpha}}d_{\mathcal{W}^{\delta}(1)}  \pa{f_J, f}+ V(J_0)}+\frac{C}{n\alpha^2}.
\end{equation}

As the term proportional to $(n\alpha^2)^{-1}$ is negligible compared to the other one, determining a $J_0$ in $\mathcal{M}_n^\alpha$ that realizes the compromise between the bias and the variance terms leads to the rates below.

\begin{corollary}\label{cor:adaptupper}
	Assume that $f$ belongs to $\mathcal{W}^{\beta}(R)$, with positive smoothness parameter $\beta$, and positive radius $R$ such that $R^2\geq d$. 
	Under the assumptions of Theorem \ref{thm:oraclebound}, the adaptive estimator $\hat{f}_{ \hat{J}}$ satisfies
	$$
	\sup_{f\in\mathcal{W}^{\beta}(R)} \esp{d_{\mathcal{W}^{\delta}(1)}\pa{\hat{f}_{ \hat{J}},f}} \leq
	\left\{ \begin{array}{ll}
		C \left(\frac{n\alpha^2}{\log(n\alpha^2)}\right)^{-\frac{\beta+\delta}{2\beta+2d}} & \text{if $\delta<d$,} \\
		C \pa{\frac{n\alpha^2}{\cro{\log(n\alpha^2)}^{4d+1}}}^{-1/2} & \text{if $\delta=d$,}\\
		C \left(\frac{n\alpha^2}{\log(n\alpha^2)}\right)^{-1/2} & \text{if $\delta>d$},
	\end{array}\right.
	$$
	where $C=C_{A,d,\delta,\beta,R}$ denotes a constant depending on $A$, $d$, $\delta$, $\beta$, and $R$. 
\end{corollary}

The proof of Corollary \ref{cor:adaptupper} is detailed in Section \ref{sect:proof:cor:adaptupper}. In particular, when compared to the lower bounds obtained in Theorem \ref{thm:lowerbound}, Corollary \ref{cor:adaptupper} ensures that $\hat{f}_{\hat{J}}$ is {\it optimal} in the minimax sense up to logarithm factors and thus adaptive. 
However, it requires a huge amount of computations. Indeed, for all $J$ in the collection $\mathcal{M}_n^\alpha$, we need to privatize again the data since the privacy mechanism $Q^{(J)}$ depends on $J$. This leads to $\sum_{J\in\mathcal{M}_n^\alpha} J^d$ coefficients to privatize for each individual $1\leq i\leq n$. \\

Note that in practice, constants $0<\kappa_1\leq\kappa_2$ are used to define the criterion, namely 
$$
\hat{A}(J)=\max_{ J'\in \mathcal{M}_n^\alpha}\left\{  \left( d_{\mathcal{W}^{\delta}(1)}  
\pa{\hat f_{J'},\hat f_{J'\wedge J}}-\kappa_1V(J')\right)_{+}
\right\},
\quad \text{and} \quad 
\widehat{\operatorname{Crit}}(J) = \hat{A}(J) + \kappa_2V(J).
$$ These constants need to be calibrated numerically through a preliminary study. One may choose $\kappa_1 = \kappa_2$ in order to reduce the number of parameters to tune, which is how the procedure was presented. However, recent work by \cite{lacour2016minimal} in kernel density estimation framework, suggests to first calibrate $\kappa_1$ and set $\kappa_2=2\kappa_1$ as $\hat{A}(J)$ and $V(J)$ exhibit antagonistic behaviours.
For simplicity, we choose to fix them equal to 2 in this study.

\section{Results over anisotropic Sobolev balls}\label{sect:anisotrop}

Assume now that the regularity is not the same in each direction. 
More precisely, consider the multidimensional parameter $\mi{\beta} = (\beta_1,\ldots,\beta_d)$ in $(\R_+^*)^d$ and define the \emph{anisotropic} Sobolev ball 

\begin{equation}\label{eq:defaniSobolev}
	\mathcal{W}^{\mi{\beta}}(R) = \ac{f\in\L_2([0,1]^d) \; ; \; 
		\sum_{\mi{j}\in(\N^*)^d} \pa{j_1^{2\beta_1} + \ldots + j_d^{2\beta_d}}\theta_{\mi{j}}^2(f)\leq R^2}.
\end{equation}

All the results may be easily adapted as soon as one considers discriminators with the same anisotropy as the density regularity, that is anisotropic Sobolev balls $\mathcal{W}_{\mi{\delta}}(1)$, where $\mi{\delta} = (\delta_1,\ldots,\delta_d)$ is such that the ratio $\beta_m/\delta_m$ is constant. 

More precisely, denote 
\begin{equation}
	\label{eq:defbetadeltaani}
	\frac{1}{\beta} = \frac{1}{d}\sum_{m=1}^d \frac{1}{\beta_m} 
	\quad \text{and}\quad
	\frac{1}{\delta} = \frac{1}{d}\sum_{m=1}^d \frac{1}{\delta_m}. 
\end{equation}
Then, $\mathcal{W}_{\mi{\beta}}(R)$ and $\mathcal{W}_{\mi{\delta}}(1)$ have the same anisotropy if and only if for all $1\leq m\leq d$, 
\begin{equation}\label{eq:sameani}
	\frac{\beta}{\beta_m} = \frac{\delta}{\delta_m}.
\end{equation}

In the anisotropic case, the choice of the number of coordinates considered in each dimension should depend on the corresponding regularity. 
Hence, for an integer $J$ to be calibrated later, denote the set of multi-indices that are considered
$$
\mathcal{J} = \prod_{m=1}^d \ac{1,\ldots,J_m}
\quad \text{where} \quad 
\forall\ 1\leq m\leq d,\ J_m \asymp J^{\beta/\beta_m} = J^{\delta/\delta_m}.
$$
In particular, the number of elements in $\mathcal{J}$ equals 
$$\# \mathcal{J} = \prod_{m=1}^d J_m = J^d.$$
Note that in the isotropic case, for all $1\leq m\leq d$,  $\beta_m=\beta$ thus $J_m = J$. Then, one recovers the set of multi-indices $\mathcal{J} = \{1,2,\ldots,J\}^d$ over which we sum in the definition of the estimator \eqref{eq:privateprojestim}.

\begin{theorem}
	\label{thm:lowerboundanisotrop}
	Consider a sample size $n$, a privacy level $\alpha$ in $(0,A]$ such that $n\alpha^2\geq 1$, integer smoothness parameters $\mi{\beta}\in(\N^*)^d$ and $\mi{\delta}\in(\N^*)^d$ with same anisotropy as defined in Equation \eqref{eq:sameani}, and a positive radius $R$ such that $R^2>d$. 
	Then, there exists a positive constant $c=c_{\mi{\beta},\mi{\delta},R,d,A}$ such that 
	$$
	\rho_n(\mathcal{W}^{\mi{\beta}}(R),\mathcal{W}^{\mi{\delta}}(1),\mathcal{Q}_\alpha) \geq
	c\max\ac{(n\alpha^2)^{-\frac{\beta+\delta}{2\beta+2d}},(n\alpha^2)^{-1/2}},
	$$
	where $\beta$ and $\delta$ are defined by Equation \eqref{eq:defbetadeltaani}. 
\end{theorem}

\paragraph{}
For the upper bound, we also need to adapt the privacy mechanism. The only difference comes from the fact that we do not consider the same number of dyadic blocks in each dimension. 
In this case, define 
$$L_m = \lfloor \log_2\pa{J^{\beta/\beta_m}+1} \rfloor -1,
\quad \text{and}\quad 
J_m = 2^{L_m+1}-1,$$ 
such that $J^{\beta/\beta_m} \leq J_m \leq 3J^{\beta/\beta_m}$. 
We then consider the same dyadic blocks $\mathcal{J}_{\mi{\ell}}$ for $\mi{\ell}$ in $\mathcal{L} = \prod_{m=1}^d \ac{0,1,\ldots,L_m}$,
as in the isotropic case, such that 
$$\mathcal{J} = \bigcup_{\mi{\ell}\in\mathcal{L}} \mathcal{J}_{\mi{\ell}}.$$ 
Then, we apply the non-interactive Coordinate dyadic block privacy mechanism described in Section \ref{sect:priv_machanism} and denote it $Q^{(J)}$. 
Finally, we define the 
private projection estimator by  
\begin{equation}
	\label{eq:privateprojestimanisotrop}
	\hat{f}_J = \sum_{\mi{j}\in\mathcal{J}}\hat\theta_{\mi{j}} \varphi_{\mi{j}},
\end{equation}
for a well calibrated $J$. 

\begin{theorem}\label{thm:upperanisotrop}
	Consider a privacy level $\alpha$ in $(0,A]$, smoothness parameters $\mi{\beta}$ in $(\R_+^*)^d$ and $\mi{\delta}$ in $(\R_+^*)^d$ with same anisotropy as defined in Equation \eqref{eq:sameani}, and a positive radius $R$ such that $R^2\geq d$. 
	Consider the Coordinate dyadic block privacy mechanism $Q^{(J)}$, where 
	$$
	\alpha_{\mi{\ell}} = \frac{\alpha}{S_{d,\delta}(J)} \prod_{m=1}^d 2^{\ell_m(1-\delta_m/d)/2} 
	\qquad \text{and}\qquad 
	S_{d,\delta}(J) = \sum_{\mi{\ell'}\in\mathcal{L}} \prod_{m=1}^d 2^{\ell_m'(1-\delta_m/d)/2}.
	$$ 
	Then, the projection estimator $\hat{f}_J$ defined in Equation \eqref{eq:privateprojestimanisotrop} satisfies 
	$$
	\sup_{f\in\mathcal{W}^{\mi{\beta}}(R)} \hspace{-10pt} \esp{d_{\mathcal{W}^{\mi{\delta}}(1)}\pa{\hat{f}_J,f}} \leq
	\left\{ \begin{array}{ll}
		C (n\alpha^2)^{-\frac{\beta+\delta}{2\beta+2d}} & \text{if $\forall\ m,\ \delta_m<d$ and }J\asymp (n\alpha^2)^{\frac{1}{2\beta+2d}} \\
		C \pa{\frac{n\alpha^2}{\log(n\alpha^2)^{4d}}}^{-1/2} \hspace{-5pt} & \text{if $\forall\ m,\ \delta_m=d$ and }J\asymp\pa{\frac{n\alpha^2}{\log(n\alpha^2)^{4d}}}^{\frac{1}{2\beta+2\delta}}\\
		C (n\alpha^2)^{-1/2} & \text{if $\forall\ m,\ \delta_m>d$ and }J\asymp(n\alpha^2)^{\frac{1}{2\beta+2\delta}},
	\end{array}\right.
	$$
	where $\beta$ and $\delta$ are defined in Equation \eqref{eq:defbetadeltaani}, and where $C=C_{A,d,\mi{\delta},\mi{\beta},R}$ denotes a constant depending on $A$, $d$, $\mi{\delta}$, $\mi{\beta}$, and $R$. 
\end{theorem}

One may notice that the upper bounds are obtained only for all $\delta_m$ on the same regime (that is all less than $d$, all equal to $d$, or all greater than $d$). This condition is due to technical reasons, and appears naturally in the proof. Without this condition, the minimax optimality remains an open question. 
Moreover, as in the isotropic case, the private projection estimator with the Coordinate dyadic block privacy mechanism are minimax optimal for all $\delta_m$ less than $d$, or all $\delta_m$ greater than $d$.

\paragraph{}
The adaptation of the proofs of Theorems \ref{thm:lowerbound} and \ref{thm:upper} to the anisotropic case with the adjustments defined above is straightforward, and leads to Theorems \ref{thm:lowerboundanisotrop} and \ref{thm:upperanisotrop}. 
For the sake of clarity, we decide to present the proofs only in the isotropic case. 
Moreover, as the choice of the regularity parameter $\mi{\delta}$ of the generators, and of $J$ which determines the number of considered coefficients both depend on the unknown regularity $\mi{\beta}$ of the density, this procedure is not adaptive in the minimax sense. This is a topic for future work.

\section{Main proofs}
\label{sect:proofs}

\subsection{Proof of Theorem \ref{thm:lowerbound}}
\label{sect:proofthmlowerbound}

Let $n\geq 1$ and $\alpha$ in $(0,A]$ such that $n\alpha^2\geq 1$. 
Let $\beta$ and $\delta$ be two positive integer smoothness parameters, and consider a positive radius such that $R^2>d$. 
Let $J\geq 1 $ be a positive integer. 
We define $\psi$ on $[0,1]$ by 
$$\psi(t) = \exp\pa{\frac{-1}{1-(4t-1)^2}}\1{t\in(0,1/2)} - \exp\pa{\frac{-1}{1-(4t-3)^2}}\1{t\in(1/2,1)}.$$
Note that the function $\psi$ is based on bump functions and often used to prove lower bounds in nonparametric statistics. 
In particular, $\psi$ is a periodic function on $[0,1]$ with continuous derivatives of all orders, such that all derivatives are uniformly bounded on $[0,1]$ and periodic, and $\psi$ satisfies $\int_{[0,1]}\psi(x)dx=0$. \\
Denote for all $\mi{j}=(j_1,\ldots,j_d)$ in $\ac{1,2,\ldots,J}^d$, and all $x=(x_1,\ldots,x_d)$ in $[0,1]^d$, 
$$
G_{\mi{j}}(x) = \prod_{m=1}^d \psi\pa{J\pa{x_m - \frac{j_m-1}{J}}}. 
$$
Note that the support of $G_{\mi{j}}$ is $\prod_{m=1}^d \cro{\frac{j_m-1}{J},\frac{j_m}{J}}$, and for all $p>0$, 
\begin{equation}\label{eq:momentsGj}
	\int_{[0,1]^d}G_{\mi{j}} = 0, 
	\qquad
	\qquad \text{and}\qquad
	\norm{G_{\mi{j}}}_{p}^p = \frac{\pa{\norm{\psi}_p^{p}}^d}{J^d}.
\end{equation}
Indeed, 
\begin{eqnarray*}
	\norm{G_{\mi{j}}}_{p}^p 
	\ =\  \int_{[0,1]^d}\abs{G_{\mi{j}}(x)}^pdx 
	&=& \prod_{m=1}^d \cro{\int_{\frac{j_m-1}{J}}^{\frac{j_m}{J}} \abs{\psi\pa{J \pa{x_m-\frac{j_m-1}{J}}}}^p dx_m} \\
	&=& \prod_{m=1}^d \cro{\frac{1}{J}\int_{0}^{1} \abs{\psi\pa{y_m}}^p dy_m} 
	\ = \ \frac{(\norm{\psi}_p^{p})^d}{J^{d}}. 
\end{eqnarray*}
Consider $\gamma_n,\eta>0$ that will be calibrated at the end of the proof, and define two parametrized families of functions on $[0,1]^d$ by
$$
\mathcal{F}^{\beta}(\gamma_n) = \ac{f_{\nu} = \1{[0,1]^d} + \frac{\gamma_n}{J^{\beta}} \sum_{\mi{j}\in\ac{1,2,\ldots,J}^d} \nu_{\mi{j}} G_{\mi{j}}; \ \nu\in\ac{0,1}^{J^d}},
$$
and
$$
\mathcal{D}^{\delta}(\eta) = \ac{g_{\lambda} = \frac{\eta}{J^\delta}\sum_{\mi{j}\in\ac{1,2,\ldots,J}^d} \lambda_{\mi{j}} G_{\mi{j}}; \ \lambda\in\ac{-1,1}^{J^d}}.
$$

\begin{lemma}\label{lm:inclusionParamdansSobolev} Let $\beta$ and $\delta$ be two positive integer parameters. Assume $R^2 > d$. 
	\begin{enumerate}
		\item All functions $f_\nu$ in $\mathcal{F}^{\beta}(\gamma_n)$ are densities as soon as 
		$\gamma_n \leq \norm{\psi}_\infty^{-d}$. 
		\item One has the inclusions  
		$\mathcal{F}^{\beta}(\gamma_n)\subset \mathcal{W}^\beta(R)$ 
		and 
		$\mathcal{D}^{\delta}(\eta)\subset \mathcal{W}^\delta(1),$
		as soon as 
		\begin{equation}\label{eq:condgammeta}
			\gamma_n^2\leq \frac{R^2 - d}{d\norm{\psi}_2^{2(d-1)}  \cro{ \norm{\psi}_2^{2} + \norm{\psi^{(\beta)}}_2^2 }}, 
			\quad \text{and}\quad
			\eta^2\leq \frac{1}{d\norm{\psi}_2^{2(d-1)}  \cro{ \norm{\psi}_2^{2} + \norm{\psi^{(\delta)}}_2^2 }}. 
		\end{equation}
	\end{enumerate}
\end{lemma}
The proof of Lemma \ref{lm:inclusionParamdansSobolev} is detailed in Section \ref{sect:proofinclusionParamdansSobolev}. \\

In the following, consider $\gamma_n$ and $\eta$ such that $\gamma_n\leq\norm{\psi}_\infty^{-d},$ and that $\eqref{eq:condgammeta}$ is satisfied.
In particular, for all $\alpha$-LDP privacy mechanisms $Q$, that satisfies Equation \eqref{eq:localPrivacy}, and all estimators $\tilde{f}$ of $f$ based on data privatized from $Q$, 
\begin{eqnarray*}
	\sup_{f\in \mathcal{W}^{\beta}(R)}\esps{f,Q}{d_{\mathcal{W}^{\delta}(1)}\pa{\tilde{f},f}} 
	&\geq & \sup_{f_\nu\in \mathcal{F}^{\beta}(\gamma_n)}\esps{f_\nu,Q}{d_{\mathcal{W}^{\delta}(1)}\pa{\tilde{f},f_\nu}} \\
	&=& \max_{\nu\in\{0,1\}^{J^d}} \esps{f_\nu,Q}{d_{\mathcal{W}^{\delta}(1)}\pa{\tilde{f},f_\nu}}.
\end{eqnarray*}
Moreover, for all $\nu$ in $\{0,1\}^{J^d}$,   
\begin{eqnarray}
	d_{\mathcal{W}^{\delta}(1)}\pa{\tilde{f},f_\nu} 
	&=& \sup_{g\in \mathcal{W}^{\delta}(1)} \int_{[0,1]^d}\cro{\pa{\tilde{f}-f_\nu}g} \nonumber\\
	&\geq & \sup_{g_\lambda\in \mathcal{D}^{\delta}(\eta)} \int_{[0,1]^d}\cro{\pa{\tilde{f}-f_\nu}g_\lambda} \nonumber\\
	&=& \max_{\lambda\in\{-1,1\}^{J^d}} \frac{\eta}{J^{\delta}} \sum_{\mi{j}\in\{1,\ldots,J\}^d}\lambda_{\mi{j}} \int_{[0,1]^d} \cro{\pa{\tilde{f}-f_\nu}G_{\mi{j}}} \nonumber\\
	&\geq & \frac{\eta}{J^{\delta}} \sum_{\mi{j}\in\{1,\ldots,J\}^d}\abs{\int_{[0,1]^d} \cro{\pa{\tilde{f}-f_\nu}G_{\mi{j}}}},\label{eq:minordadvser}
\end{eqnarray}
with the particular choice of $\tilde\lambda$ such that for all $\mi{j}$, $\tilde\lambda_{\mi{j}}$ is the sign of $\int_{[0,1]^d} [(\tilde{f}-f_\nu)G_{\mi{j}}]$. 
Yet, by definition of $f_\nu$ in $\mathcal{F}^\beta(\gamma_n)$, and since the supports of the functions $G_{\mi{j}}$ are disjoint, 
$$
\abs{\int_{[0,1]^d} \cro{\pa{\tilde{f}-f_\nu}G_{\mi{j}}}}
= \abs{\int_{[0,1]^d} \cro{\pa{\tilde{f}-h_{\mi{j}}(\nu_{\mi{j}})}G_{\mi{j}}}},
$$
where $h_{\mi{j}}(\nu_{\mi{j}}) = 1+\frac{\gamma_n}{J^{\beta}} \nu_{\mi{j}} G_{\mi{j}}$. 
Introduce $$\tilde{\nu}_{\mi{j}} \in \argmin_{\nu_{\mi{j}}\in\{0,1\}} \abs{\int_{[0,1]^d} \cro{\pa{\tilde{f}-h_{\mi{j}}(\nu_{\mi{j}})}G_{\mi{j}}}}.$$
Then by the triangular inequality, 
\begin{eqnarray}
	\abs{\int_{[0,1]^d} \cro{\pa{\tilde{f}-f_\nu}G_{\mi{j}}}}
	&\geq& \frac{1}{2}\ac{\abs{\int_{[0,1]^d} \cro{\pa{\tilde{f}-h_{\mi{j}}(\tilde{\nu}_{\mi{j}})}G_{\mi{j}}}}+\abs{\int_{[0,1]^d} \cro{\pa{\tilde{f}-h_{\mi{j}}(\nu_{\mi{j}})}G_{\mi{j}}}}} \nonumber\\
	&\geq& \frac{1}{2}\abs{\int_{[0,1]^d} \cro{\pa{h_{\mi{j}}(\tilde{\nu}_{\mi{j}})-h_{\mi{j}}(\nu_{\mi{j}})}G_{\mi{j}}}} \nonumber\\
	&=& \frac{\gamma_n}{2J^\beta}\abs{\tilde{\nu}_{\mi{j}} - \nu_{\mi{j}}} \int_{[0,1]^d} G^2_{\mi{j}} \nonumber\\
	&=& \frac{\norm{\psi}_2^{2d}}{2}\frac{\gamma_n}{J^{\beta+d}}\abs{\tilde{\nu}_{\mi{j}} - \nu_{\mi{j}}}, \label{eq:minorunj}
\end{eqnarray}
by Equation \eqref{eq:momentsGj} for $p=2$. 
Finally, combining \eqref{eq:minordadvser} and \eqref{eq:minorunj} leads to 
$$
d_{\mathcal{W}^{\delta}(1)}\pa{\tilde{f},f_\nu} 
\geq \frac{\eta\norm{\psi}_2^{2d}}{2}\frac{\gamma_n}{J^{\beta+\delta+d}} \rho_{H}(\tilde{\nu},\nu),
$$
where $\rho_{H}(\tilde{\nu},\nu) = \sum_{\mi{j}}|\tilde{\nu}_{\mi{j}}-\nu_{\mi{j}}| = \sum_{\mi{j}}\1{\tilde{\nu}_{\mi{j}} \neq \nu_{\mi{j}}}$ denotes the Hamming distance between $\tilde{\nu}$ and $\nu$ in $\{0,1\}^{J^d}$. Thus,
\begin{eqnarray}
	\sup_{f\in \mathcal{W}^{\beta}(R)}\esps{f,Q}{d_{\mathcal{W}^{\delta}(1)}\pa{\tilde{f},f}} 
	&\geq & \frac{\eta\norm{\psi}_2^{2d}}{2}\frac{\gamma_n}{J^{\beta+\delta+d}} \max_{\nu\in \{0,1\}^{J^d}}\esps{f_\nu,Q}{\rho_{H}(\tilde{\nu},\nu)} \nonumber\\
	&\geq & \frac{\eta\norm{\psi}_2^{2d}}{2}\frac{\gamma_n}{J^{\beta+\delta+d}} \times \inf_{\hat\nu}\max_{\nu\in \{0,1\}^{J^d}}\esps{f_\nu,Q}{\rho_{H}(\hat{\nu},\nu)}, \label{eq:ineqgeneralavanttsybakov}
\end{eqnarray}
where the infimum is taken over all estimators $\hat\nu$ based on data from the distribution $P_\nu$ with density $f_\nu$, that have been privatized  using a privacy mechanism $Q$ in $\mathcal{Q}_\alpha$. 
Denote $M^n_{\nu}$ the distribution of such privatized data. 
In order to lower bound this infimum, let us apply Theorem 2.12 (p. 118) of \cite{tsybakov2009introduction} recalled below. 

\begin{theorem}[\cite{tsybakov2009introduction}, Theorem 2.12] \label{lm:tsybakovthm2.12}
	Let $\Theta = \{0,1\}^N$, with $N\geq 1$ and $\ac{P_\theta,\theta\in\Theta}$ be a set of $2^N$ probability measures on a measurable space $(\mathcal{X},\mathcal{A})$. Denote $\mathds{E}_\theta$ the corresponding expectations. 
	Assume that there exists $\xi>0$ such that for all $\theta,\theta'$ in $\Theta$ satisfying $\rho_{H}(\theta,\theta')=1$, the Kullback-Leibler divergence between $P_\theta$ and $P_{\theta'}$ satisfies $KL(P_\theta,P_{\theta'}) \leq \xi$. 
	Then, it yields 
	$$\inf_{\hat\theta}\max_{\theta\in\Theta}\esps{\theta}{\rho_{H}\pa{\hat\theta,\theta}} \geq \frac{N}{2} \max\ac{\frac{e^{-\xi}}{2},1-\sqrt{\frac{\xi}{2}}}.$$ 
\end{theorem}

We thus need to upper bound the Kullback-Leibler divergence of the privatized distributions $M^n_{\nu}$ and $M^n_{\nu'}$, where $\nu$ and $\nu'$ belong to $\{0,1\}^{J^d}$ such that $\rho_{H}(\nu,\nu')=1$. 
Yet, according to Corollary 3 of \cite{duchi2018minimax} applied in the particular case of i.i.d. $X_i$'s, 
$$\KL(M^n_{\nu},M^n_{\nu'}) \leq 4 n \pa{e^\alpha-1}^2 \TV^2(P_\nu,P_{\nu'}).$$
Moreover, if $\rho_{H}(\nu,\nu')=1$, then there exists a unique $\mi{j_0}$ such that $\nu_{\mi{j_0}} \neq \nu_{\mi{j_0}}'$, and thus
$$
\TV(P_\nu,P_{\nu'}) = \frac{1}{2}\int_{[0,1]^d}\abs{f_\nu-f_{\nu'}} 
= \frac{\gamma_n}{2J^\beta}\underbrace{\abs{\nu_{\mi{j_0}} - \nu_{\mi{j_0}}'}}_{=1}\int_{[0,1]^d}\abs{G_{\mi{j_0}}} 
= \frac{\norm{\psi}_1^{d}}{2}\frac{\gamma_n}{J^{\beta+d}},
$$
by Equation \eqref{eq:momentsGj} for $p=1$. 
We deduce that 
$$\KL(M^n_{\nu},M^n_{\nu'}) \leq \norm{\psi}_1^{2d}\frac{\gamma_n^2}{J^{2\beta+2d}}n \pa{e^\alpha-1}^2
\leq e^{2A}\norm{\psi}_1^{2d} \frac{\gamma_n^2}{J^{2\beta+2d}}n\alpha^2,$$
as one can prove by Taylor Theorem with the Lagrange form of remainder that for all $\alpha$ in $(0,A]$, $e^\alpha-1 \leq \alpha e^A.$
In particular, if 
$ \frac{\gamma_n^2}{J^{2\beta+2d}}n\alpha^2 =: \gamma^2$
is constant, we obtain that 
$$\KL(M^n_{\nu},M^n_{\nu'}) \leq \xi, \qquad \text{where}\qquad \xi = \gamma^2 e^{2A}\norm{\psi}_1^{2d}.$$
Applying Theorem \ref{lm:tsybakovthm2.12} with $N=J^{d}$ to the privatized distributions leads to 
$$\inf_{\hat\nu}\max_{\nu\in\{0,1\}^{J^d}}\esps{f_\nu,Q}{\rho_{H}\pa{\hat\nu,\nu}} \geq \frac{J^d}{2} \max\ac{\frac{e^{-\xi}}{2},1-\sqrt{\frac{\xi}{2}}}.$$ 
Finally, from \eqref{eq:ineqgeneralavanttsybakov}, we obtain 
$$
\sup_{f\in \mathcal{W}^{\beta}(R)}\esps{f,Q}{d_{\mathcal{W}^{\delta}(1)}\pa{\tilde{f},f}} 
\geq c_0 \frac{\gamma_n}{J^{\beta+\delta}},
$$
where $c_0 = \pa{\eta\norm{\psi}_2^{2d}/4}\max\ac{e^{-\xi}/2, 1-\sqrt{\xi/2}}$. 
This being true for all $Q$ in $\mathcal{Q}_\alpha$ and all $\tilde{f}$, we deduce that 
$$\rho_n(\mathcal{W}^{\beta}(R),\mathcal{W}^{\delta}(1),\mathcal{Q}_\alpha) = \inf_{Q\in\mathcal{Q}_\alpha}\inf_{\tilde{f}}\sup_{f\in\mathcal{W}^{\beta}(R)} \esp{d_{\mathcal{W}^{\delta}(1)}\pa{\tilde{f},f}} \geq c_0 \frac{\gamma_n}{J^{\beta+\delta}}.
$$
In order to obtain both rates in Theorem \ref{thm:lowerbound}, consider $\gamma$ and $\eta$ such that $$\gamma^2 = \min\ac{\norm{\psi}_\infty^{-2d}, 
	\frac{R^2 - d}{d\norm{\psi}_2^{2(d-1)}  \cro{ \norm{\psi}_2^{2} + \norm{\psi^{(\beta)}}_2^2 }}
},
\quad\text{and}\quad 
\eta^2 = \frac{1}{d\norm{\psi}_2^{2(d-1)}  \cro{ \norm{\psi}_2^{2} + \norm{\psi^{(\delta)}}_2^2 }} $$
so that we may apply Lemma \ref{lm:inclusionParamdansSobolev} for the different choices of $\gamma_n$ below. 
\begin{itemize}
	\item First, choosing $J = (n\alpha^2)^{\frac{1}{2\beta+2d}}$ and $\gamma_n = \gamma$ implies that $\gamma_n^2 \frac{n\alpha^2}{J^{2\beta+2d}} = \gamma^2$ is a constant and thus, 
	$$
	\rho_n(\mathcal{W}^{\beta}(R),\mathcal{W}^{\delta}(1),\mathcal{Q}_\alpha) 
	\geq c_0\frac{\gamma}{J^{\beta+\delta}} 
	= c\ (n\alpha^2)^{-\frac{\beta+\delta}{2\beta+2d}}.
	$$
	\item Second, choosing $J=1$ and $\gamma_n = \gamma/(\sqrt{n\alpha^2})\leq \gamma$ (as $n\alpha^2\geq1$) also implies that 
	$\gamma_n^2 \frac{n\alpha^2}{J^{2\beta+2d}} = \gamma^2$ is a constant and thus, 
	$$
	\rho_n(\mathcal{W}^{\beta}(R),\mathcal{W}^{\delta}(1),\mathcal{Q}_\alpha) 
	\geq c_0\frac{\gamma}{\sqrt{n\alpha^2}} 
	= c\ (n\alpha^2)^{-1/2},
	$$
	where $c$ is a constant depending on $\beta,\delta,R,d,A$. \end{itemize}
This concludes the proof of Theorem \ref{thm:lowerbound}. 

\subsection{Proof of Proposition \ref{prop:CBmechanismalphaLDP}}
\label{sect:proof:prop:CBmechanismalphaLDP}

Before proving Propositions \ref{prop:CBmechanismalphaLDP} and \ref{prop:CBmechanismespvar}, let us introduce some useful notation. 
Define for all $\ml$ in $\mcL$ and for all $\tilde{v}$ in $\{-B_0,B_0\}^{d_{\ml}}$, 
\begin{eqnarray*}
	\mathcal{D}_{\ml}^+\pa{\tilde{v}} &=& \ac{z\in\{\pm B_{d_{\mi{\ell}}}(\alpha_{\mi{\ell}})\}^{d_{\mi{\ell}}} \; ; \; 
		\langle z,\tilde{v} \rangle >0},\\
	\mathcal{D}_{\ml}^-\pa{\tilde{v}} &=& \ac{z\in\{\pm B_{d_{\mi{\ell}}}(\alpha_{\mi{\ell}})\}^{d_{\mi{\ell}}} \; ; \; 
		\langle z,\tilde{v} \rangle <0},\\
	\mathcal{D}_{\ml}^0\pa{\tilde{v}} &=& \ac{z\in\{\pm B_{d_{\mi{\ell}}}(\alpha_{\mi{\ell}})\}^{d_{\mi{\ell}}} \; ; \; 
		\langle z,\tilde{v} \rangle =0},
\end{eqnarray*}
where $B_{d_{\mi{\ell}}}(\alpha_{\mi{\ell}})$ is defined in Equation \eqref{eq:defBdalpha}.
In particular, the sets $\mathcal{D}_{\ml}^+\pa{\tilde{v}}$, $\mathcal{D}_{\ml}^-\pa{\tilde{v}}$, and $\mathcal{D}_{\ml}^0\pa{\tilde{v}}$ form a partition of $\{\pm B_{d_{\mi{\ell}}}(\alpha_{\mi{\ell}})\}^{d_{\mi{\ell}}}$. 
Moreover, 
\begin{equation}
	\label{eq:lienD+D-} 
	z \in \mathcal{D}_{\ml}^+\pa{\tilde{v}} \quad \Longleftrightarrow (-z) \in \mathcal{D}_{\ml}^-\pa{\tilde{v}}.
\end{equation}
We deduce that $\abs{\mathcal{D}_{\ml}^+\pa{\tilde{v}}} = \abs{\mathcal{D}_{\ml}^-\pa{\tilde{v}}}$ and thus 
$$2\abs{\mathcal{D}_{\ml}^+\pa{\tilde{v}}} + \abs{\mathcal{D}_{\ml}^0\pa{\tilde{v}}} = 2^{d_{\ml}}.$$
Besides, if $d_{\ml}$ is even, $\abs{\mathcal{D}_{\ml}^0\pa{\tilde{v}}} = \binom{d_{\ml}}{d_{\ml}/2}$ since $z$ belongs to $\mathcal{D}_{\ml}^0\pa{\tilde{v}}$ if and only if half of its coordinates have the same sign as $\tilde{v}$, and the other half have the opposite sign. 
On the contrary, if $d_{\ml}$ is odd, then $\abs{\mathcal{D}_{\ml}^+\pa{\tilde{v}}} = 2^{d_{\ml}-1}$ since for all $\tilde{v}$ in $\{\pm B_0\}^{d_{\ml}}$, $\mathcal{D}_{\ml}^0\pa{\tilde{v}} = \varnothing$. 
Recalling the definition of $p_{\ml}$, we deduce that in all cases, 
\begin{equation}\label{eq:calculpell}
	p_{\ml} 
	= 1 - \frac{\abs{\mathcal{D}_{\ml}^0\pa{\tilde{v}}}}{2^{d_{\ml}}}
	= \frac{\abs{\mathcal{D}_{\ml}^+\pa{\tilde{v}}}}{2^{d_{\ml}-1}}. 
\end{equation}

To prove the $\alpha$-local differential privacy, it is sufficient to prove the $\alpha$-local differential privacy in each bloc, that is for all 
$$z = \pa{z_{[\mi{\ell}]}}_{\ml\in\mcL} \in \prod_{\ml\in\mcL}\ac{\pm B_{d_{\mi{\ell}}}(\alpha_{\mi{\ell}})}^{d_{\ml}}, $$
for all $x$ and $x'$, and for all $\ml$ in $\mcL$, 
\begin{equation}\label{eq:alphablockLDP}
	\frac{\proba{Z_{[\mi{\ell}]}=z_{[\mi{\ell}]}\middle| X=x}}{\proba{Z_{[\mi{\ell}]}=z_{[\mi{\ell}]}\middle| X=x'}} \leq e^{\alpha_{\mi{\ell}}}.
\end{equation}
Indeed, by independence between the random draws $Z_{[\ml]}$ in each block, we deduce from both \eqref{eq:alphablockLDP} and $\sum_{\mi{\ell}\in\mcL}\alpha_{\mi{\ell}} = \alpha$ that 
\begin{eqnarray*}
	\frac{\proba{Z=z\middle| X=x}}{\proba{Z=z\middle| X=x'}} 
	&=& \prod_{\mi{\ell}\in\ac{0,1,\ldots,L}^d} \frac{\proba{Z_{[\mi{\ell}]}=z_{[\mi{\ell}]}\middle| X=x}}{\proba{Z_{[\mi{\ell}]}=z_{[\mi{\ell}]}\middle| X=x'}} \\
	&\leq& \prod_{\mi{\ell}\in\ac{0,1,\ldots,L}^d} e^{\alpha_{\mi{\ell}}} \ = \ e^{\alpha}.
\end{eqnarray*}
This proves that the \emph{Coordinate block privacy mechanism} $Q$ satisfies non-interactive $\alpha$-local differential privacy as in \eqref{eq:localPrivacyNonInter}. \\

Let us now prove \eqref{eq:alphablockLDP}.  Fix $x$, $x'$, $z=(z_{[\ml]})_{\ml}$ and $\ml$ in $\mcL$. Then, as conditionally on $\tilde{V}_{[\ml]}$, $Z_{[\ml]}$ is independent on $X$, 
\begin{equation}\label{eq:probaZsachantX}
	\proba{Z_{[\ml]} = z_{[\ml]}\middle|X=x} 
	= \sum_{\tilde{v}_{[\ml]}\in \{\pm B_0\}^{d_{\ml}}} \proba{Z_{[\ml]} = z_{[\ml]}\middle|\tilde{V}_{[\ml]} = \tilde{v}_{[\ml]}} \proba{\tilde{V}_{[\ml]} = \tilde{v}_{[\ml]}\middle|X=x}.
\end{equation}
We thus need to compute $\proba{Z_{[\ml]} = z_{[\ml]}\middle|\tilde{V}_{[\ml]} = \tilde{v}_{[\ml]}}$. More precisely, let us prove that 
\begin{equation}
	\label{eq:probaZsacantvtilde}
	\proba{Z_{[\ml]} = z_{[\ml]}\middle|\tilde{V}_{[\ml]} = \tilde{v}_{[\ml]}} 
	= \frac{\pi_{\alpha_{\ml}}}{2^{d_{\ml}-1}} \1{\mathcal{D}_{\ml}^+\pa{\tilde{v}_{[\ml]}}}(z_{[\ml]})
	+ \frac{1-\pi_{\alpha_{\ml}}}{2^{d_{\ml}-1}} \1{\mathcal{D}_{\ml}^-\pa{\tilde{v}_{[\ml]}}}(z_{[\ml]})
	+ \frac{1}{2^{d_{\ml}}} \1{\mathcal{D}_{\ml}^0\pa{\tilde{v}_{[\ml]}}}(z_{[\ml]}).
\end{equation}
To do this, we need to distinguish between two cases, depending on the parity of $d_{\ml}$.
\begin{itemize}
	\item If $d_{\ml}$ is even, then, since $Y_{\ml} \sim\mathcal{B}(p_{\ml})$ with $0<p_{\ml}<1$, 
	\begin{eqnarray}
		\proba{Z_{[\ml]} = z_{[\ml]}\middle|\tilde{V}_{[\ml]} = \tilde{v}_{[\ml]}} 
		& = & \ p_{\ml}\proba{Z_{[\ml]} = z_{[\ml]}\middle|\tilde{V}_{[\ml]} = \tilde{v}_{[\ml]},Y_{\ml}=1} \nonumber\\
		&& +\ (1-p_{\ml}) \proba{Z_{[\ml]} = z_{[\ml]}\middle|\tilde{V}_{[\ml]} = \tilde{v}_{[\ml]}, Y_{\ml}=0}.\label{eq:probacondYpair}
	\end{eqnarray}
	Besides, recall that conditionally on $\{\tilde{V}_{[\ml]} = \tilde{v}_{[\ml]},Y_{\ml}=1\}$, 
	$$
	T_{\ml} \sim\mathcal{B}\pa{\pi_{\alpha_{\ml}}},
	\quad \text{and}\quad
	Z_{[\mi{\ell}]} \sim \left\{\begin{array}{ll}
		\mathcal{U}\pa{\mathcal{D}_{\ml}^+\pa{\tilde{v}}} 
		& \text{if }T_{\mi{\ell}}=1,\\
		\mathcal{U}\pa{\mathcal{D}_{\ml}^-\pa{\tilde{v}}} 
		& \text{if }T_{\mi{\ell}}=0, 
	\end{array}\right.
	$$
	with $0<\pi_{\alpha_{\ml}}<1$. 
	Hence, conditionally on $\{\tilde{V}_{[\ml]} = \tilde{v}_{[\ml]},Y_{\ml}=1,T_{\ml}=1\}$, the event $\{Z_{[\ml]} = z_{[\ml]}\}$ is possible only if $z_{[\ml]}$ belongs to $\mathcal{D}_{\ml}^+\pa{\tilde{v}}$,
	and similarly when conditioning on $\{T_{\ml}=0\}$. We deduce that 
	\begin{equation}\label{eq:probacondY1pair}
		\proba{Z_{[\ml]} = z_{[\ml]}\middle|\tilde{V}_{[\ml]} = \tilde{v}_{[\ml]},Y_{\ml}=1} 
		= \frac{\pi_{\alpha_{\ml}}}{\abs{\mathcal{D}_{\ml}^+\pa{\tilde{v}_{[\ml]}}}}\1{\mathcal{D}_{\ml}^+\pa{\tilde{v}_{[\ml]}}}(z_{[\ml]})
		+ \frac{1-\pi_{\alpha_{\ml}}}{\abs{\mathcal{D}_{\ml}^-\pa{\tilde{v}_{[\ml]}}}}\1{\mathcal{D}_{\ml}^-\pa{\tilde{v}_{[\ml]}}}(z_{[\ml]}).
	\end{equation}
	Analogously, recall that conditionally on $\{\tilde{V}_{[\ml]} = \tilde{v}_{[\ml]},Y_{\ml}=0\}$, 
	$Z_{[\mi{\ell}]} \sim \mathcal{U}\pa{\mathcal{D}_{\ml}^0\pa{\tilde{v}_{[\ml]}}} $.
	Thus 
	$$\proba{Z_{[\ml]} = z_{[\ml]}\middle|\tilde{V}_{[\ml]} = \tilde{v}_{[\ml]},Y_{\ml}=0} 
	= \frac{1}{\abs{\mathcal{D}_{\ml}^0\pa{\tilde{v}_{[\ml]}}}}\1{\mathcal{D}_{\ml}^0\pa{\tilde{v}_{[\ml]}}}(z_{[\ml]}).
	$$ 
	Therefore, we deduce \eqref{eq:probaZsacantvtilde} from \eqref{eq:calculpell} and \eqref{eq:probacondYpair}. 
	\item If $d_{\ml}$ is odd, recall that $\abs{\mathcal{D}_{\ml}^+\pa{\tilde{v}_{[\ml]}}} = \abs{\mathcal{D}_{\ml}^-\pa{\tilde{v}_{[\ml]}}} = 2^{d_{\ml}-1}$, and $p_{\ml}=1$. In particular, $Y_{\ml}=1$, and we deduce from \eqref{eq:probacondY1pair} that  
	\begin{eqnarray*}
		\proba{Z_{[\ml]} = z_{[\ml]}\middle|\tilde{V}_{[\ml]} = \tilde{v}_{[\ml]}} 
		&=& \proba{Z_{[\ml]} = z_{[\ml]}\middle|\tilde{V}_{[\ml]} = \tilde{v}_{[\ml]},Y_{\ml}=1} \\
		&=& \frac{\pi_{\alpha_{\ml}}}{2^{d_{\ml}-1}} \1{\mathcal{D}_{\ml}^+\pa{\tilde{v}_{[\ml]}}}(z_{[\ml]})
		+ \frac{1-\pi_{\alpha_{\ml}}}{2^{d_{\ml}-1}} \1{\mathcal{D}_{\ml}^-\pa{\tilde{v}_{[\ml]}}}(z_{[\ml]}),
	\end{eqnarray*}
	which is equivalent to \eqref{eq:probaZsacantvtilde} since $\mathcal{D}_{\ml}^0\pa{\tilde{v}_{[\ml]}}=\varnothing$ in this case. 
\end{itemize} 
Now, notice that since $1/2 < \pi_{\alpha_{\ml}} < 1$, then 
$$\frac{1-\pi_{\alpha_{\ml}}}{2^{d_{\ml}-1}}<\frac{1}{2^{d_{\ml}}}<\frac{\pi_{\alpha_{\ml}}}{2^{d_{\ml}-1}},$$
and thus, using \eqref{eq:probaZsacantvtilde},
$$\frac{1-\pi_{\alpha_{\ml}}}{2^{d_{\ml}-1}} 
\leq \proba{Z_{[\ml]} = z_{[\ml]}\middle|\tilde{V}_{[\ml]} = \tilde{v}_{[\ml]}} 
\leq \frac{\pi_{\alpha_{\ml}}}{2^{d_{\ml}-1}}.$$
Finally, as $\sum_{\tilde{v}_{[\ml]}\in \{\pm B_0\}^{d_{\ml}}} \proba{\tilde{V}_{[\ml]} = \tilde{v}_{[\ml]}\middle|X=x} = 1$ in \eqref{eq:probaZsachantX}, we proved that 
$$\frac{1-\pi_{\alpha_{\ml}}}{2^{d_{\ml}-1}} 
\leq \proba{Z_{[\ml]} = z_{[\ml]}\middle|X=x}
\leq \frac{\pi_{\alpha_{\ml}}}{2^{d_{\ml}-1}}.$$
In particular, this implies that for all $x$ and $x'$, by definition of $\pi_{\alpha_{\ml}}$, 
$$
\frac{\proba{Z_{[\mi{\ell}]}=z_{[\mi{\ell}]}\middle| X=x}}{\proba{Z_{[\mi{\ell}]}=z_{[\mi{\ell}]}\middle| X=x'}} 
\leq \frac{\pi_{\alpha_{\ml}}}{1-\pi_{\alpha_{\ml}}} 
= \frac{\quad \frac{e^{\alpha_{\mi{\ell}}}}{1+e^{\alpha_{\mi{\ell}}}}\quad }{\frac{1}{1+e^{\alpha_{\mi{\ell}}}}} 
= e^{\alpha_{\mi{\ell}}},
$$
which ends the proof of \eqref{eq:alphablockLDP}.

\subsection{Proof of Proposition \ref{prop:CBmechanismespvar}}
\label{sect:proof:prop:CBmechanismespvar}

Fix $\mi{\ell}$ in $\{0,\ldots,L\}^d$, and $\mi{j}$ in $\mathcal{J}_{\mi{\ell}}$. 
Consider the notation $\mathcal{D}_{\ml}^+\pa{\tilde{v}}$, $\mathcal{D}_{\ml}^-\pa{\tilde{v}}$ and $\mathcal{D}_{\ml}^0\pa{\tilde{v}}$ introduced in the beginning of the proof of Proposition \ref{prop:CBmechanismalphaLDP} in Section \ref{sect:proof:prop:CBmechanismalphaLDP}. 
\begin{enumerate}
	\item 
	Let us first compute $\esp{Z_{\mi{j}}\middle|X=x}.$ 
	Since conditionally on $\tilde{V}_{[\ml]}$, $Z_{[\ml]}$ is independent on $X$, 
	\begin{eqnarray*}
		\esp{Z_{[\ml]}\middle|X=x} 
		&=& \sum_{\tilde{v}_{[\ml]}\in \{\pm B_0\}^{d_{\ml}}} \esp{Z_{[\ml]}\middle|\tilde{V}_{[\ml]} = \tilde{v}_{[\ml]}} \proba{\tilde{V}_{[\ml]} = \tilde{v}_{[\ml]}\middle|X=x}.
	\end{eqnarray*}
	Moreover, from \eqref{eq:probaZsacantvtilde}, we deduce that 
	\begin{eqnarray*}
		\esp{Z_{[\ml]}\middle|\tilde{V}_{[\ml]} = \tilde{v}_{[\ml]}} 
		&=& \sum_{z_{[\ml]}\in\{\pm B_{d_{\mi{\ell}}}(\alpha_{\mi{\ell}})\}^{d_{\mi{\ell}}}} z_{[\ml]} \proba{Z_{[\ml]} = z_{[\ml]}\middle|\tilde{V}_{[\ml]} = \tilde{v}_{[\ml]}} \nonumber \\
		&=& \frac{\pi_{\alpha_{\ml}}}{2^{d_{\ml}-1}}\ S_{\ml}^+\!\pa{\tilde{v}_{[\ml]}}
		+ \frac{1-\pi_{\alpha_{\ml}}}{2^{d_{\ml}-1}}\ S_{\ml}^-\!\pa{\tilde{v}_{[\ml]}}
		+ \frac{1}{2^{d_{\ml}}}\ S_{\ml}^0\!\pa{\tilde{v}_{[\ml]}},
	\end{eqnarray*}
	where for any $\star$ in $\{+,-,0\}$, $$S_{\ml}^{\star}\!\pa{\tilde{v}_{[\ml]}} = \sum_{z_{[\ml]}\in \mathcal{D}_{\ml}^{\star}\pa{\tilde{v}_{[\ml]}}} z_{[\ml]}.$$
	Yet, $S_{\ml}^0\!\pa{\tilde{v}_{[\ml]}} = 0$ since for all $z_{[\ml]}$ in $\mathcal{D}_{\ml}^0\pa{\tilde{v}_{[\ml]}}$, $(-z_{[\ml]})$ also belongs to $\mathcal{D}_{\ml}^0\pa{\tilde{v}_{[\ml]}}$. 
	Moreover, by \eqref{eq:lienD+D-}, $S_{\ml}^-\!\pa{\tilde{v}_{[\ml]}} = -S_{\ml}^+\!\pa{\tilde{v}_{[\ml]}}$. 
	Therefore, 
	$$
	\esp{Z_{[\ml]}\middle|\tilde{V}_{[\ml]} = \tilde{v}_{[\ml]}} 
	=\frac{2\pi_{\alpha_{\ml}}-1}{2^{d_{\ml}-1}} S_{\ml}^+\!\pa{\tilde{v}_{[\ml]}}. 
	$$
	Furthermore, \cite{butucea2023phase} proved in Section B.2 that 
	\begin{equation*}
		S_{\ml}^+\!\pa{\tilde{v}_{[\ml]}}
		= \left\{\begin{array}{ll}
			\displaystyle\frac{B_{d_{\ml}}(\alpha_{\ml})}{B_0} \binom{d_{\ml}-1}{(d_{\ml}-1)/2}\tilde{v}_{[\ml]} & \text{if $d_{\ml}$ is odd}, \\ & \\
			\displaystyle\frac{B_{d_{\ml}}(\alpha_{\ml})}{B_0} \binom{d_{\ml}-1}{d_{\ml}/2-1}\tilde{v}_{[\ml]} & \text{if $d_{\ml}$ is even},\\
		\end{array}\right. 
		= \frac{B_{d_{\ml}}(\alpha_{\ml})}{B_0} \binom{d_{\ml}-1}{\lfloor \frac{d_{\ml}-1}{2}\rfloor}\tilde{v}_{[\ml]}. 
	\end{equation*}
	Besides, 
	$$
	2\pi_{\alpha_{\ml}}-1 
	= \frac{2e^{\alpha_{\mi{\ell}}}}{1+e^{\alpha_{\mi{\ell}}}} - \frac{1+e^{\alpha_{\mi{\ell}}}}{1+e^{\alpha_{\mi{\ell}}}} 
	= \frac{e^{\alpha_{\mi{\ell}}}-1}{e^{\alpha_{\mi{\ell}}}+1}.
	$$
	Therefore, by definition of $B_{d_{\ml}}(\alpha_{\ml})$ in Equation \eqref{eq:defBdalpha}, we obtain that 
	$$
	\esp{Z_{[\ml]}\middle|\tilde{V}_{[\ml]} = \tilde{v}_{[\ml]}} 
	\ =\ \frac{1}{2^{d_{\ml}-1}} 
	\pa{\frac{e^{\alpha_{\mi{\ell}}}-1}{e^{\alpha_{\mi{\ell}}}+1}} 
	\frac{B_{d_{\ml}}(\alpha_{\ml})}{B_0}
	\binom{d_{\ml}-1}{\lfloor \frac{d_{\ml}-1}{2}\rfloor}\tilde{v}_{[\ml]} 
	\ =\ \tilde{v}_{[\ml]}.
	$$
	Finally, 
	$$
	\esp{Z_{[\ml]}\middle|X=x} 
	= \sum_{\tilde{v}_{[\ml]}\in \{\pm B_0\}^{d_{\ml}}} \tilde{v}_{[\ml]} \proba{\tilde{V}_{[\ml]} = \tilde{v}_{[\ml]}\middle|X=x}
	= \esp{\tilde{V}_{[\ml]}\middle|X=x}, 
	$$
	and thus, for all $\mj$ in $\mcJ_{\ml}$, 
	\begin{eqnarray*}
		\esp{Z_{\mj}\middle|X=x} 
		&=& \esp{\tilde{V}_{\mj}\middle|X=x} \\
		&=& B_0 \pa{\frac{1}{2} + \frac{\varphi_{\mi{j}}(x)}{2B_0}} - B_0 \pa{\frac{1}{2} - \frac{\varphi_{\mi{j}}(x)}{2B_0}} 
		= \varphi_{\mi{j}}(x).
	\end{eqnarray*}

	\item 
	In order to prove that the $Z_{\mi{j}}$ is bounded, we use the technical Lemma below. 
	\begin{lemma}\label{lm:controlBka}
		Consider $a$ in $(0,A]$ and some positive integer $k$. 
		Then, $B_k(a)$ defined as in \eqref{eq:defBdalpha} satisfies 
		$$B_k(a) \leq 2 B_0\xi_A\frac{\sqrt{k}}{a},$$
		where $\xi_A=A(e^A+1)/(e^A-1)$. 
	\end{lemma}
	The proof is detailed in Section \ref{sect:proof:lm:controlBka}.
	The upper bound on $|Z_{\mi{j}}|$ is a direct consequence of Lemma \ref{lm:controlBka}. 
	We deduce that 
	$$\var{Z_{\mi{j}}} \leq \esp{Z_{\mi{j}}^2} \leq (2 B_0\xi_A)^2\frac{d_{\mi{\ell}}}{\alpha_{\mi{\ell}}^2},$$
	which ends the proof of Proposition \ref{prop:CBmechanismespvar}. 
\end{enumerate}

\subsection{Proof of Proposition \ref{prop:controlvar}}
\label{sect:proof:prop:controlvar}

In this proof, denote for all $\mi{j}=(j_1,\ldots,j_d)$ in $(\N^*)^d$, $\norm{\mi{j}^\delta}^2 = j_1^{2\delta} + j_2^{2\delta} + \ldots +j_d^{2\delta}$.  
We aim at upper bounding 
$$
\esp{d_{\mathcal{W}^{\delta}(1)}\pa{\hat{f}_J,f_J}} = \esp{\sup_{g\in\mathcal{W}^\delta(1)}\int_{[0,1]^d}(\hat{f}_J - f_J)g}.
$$
Let $g$ in $\mathcal{W}^\delta(1)$ and consider $g=\sum_{\mi{j}\in(\N^*)^d}\theta_{\mi{j}}(g) \varphi_{\mi{j}}$ its decomposition in the Fourier basis.
Then, as the blocks $\mathcal{J}_{\mi{\ell}}$ with $\mi{\ell}$ in $\{0,\ldots,L\}^d$ form a partition of $\{1,\ldots,J\}^d$, and by the Cauchy-Schwarz inequality for all $\mi{\ell}$, 
\begin{eqnarray}
	\int_{[0,1]^d} (\hat{f}_J-f_J)g 
	& = & \sum_{\mi{j}\in\{1,\ldots,J\}^d} \cro{\hat\theta_{\mi{j}} - \theta_{\mi{j}}(f)}\theta_{\mi{j}}(g) \nonumber\\
	& = &\sum_{\mi{\ell}\in \{0,1,\ldots,L\}^d}\sum_{\mi{j}\in \mathcal{J}_{\mi{\ell}}} \cro{\hat\theta_{\mi{j}} - \theta_{\mi{j}}(f)}\frac{1}{\norm{\mi{j}^\delta}}\times\norm{\mi{j}^\delta}\theta_{\mi{j}}(g) \nonumber\\
	&\leq & \sum_{\mi{\ell}\in \{0,1,\ldots,L\}^d}\sqrt{\pa{\sum_{\mi{j}\in\mathcal{J}_{\mi{\ell}}}\cro{\hat\theta_{\mi{j}} - \theta_{\mi{j}}(f)}^2\frac{1}{\norm{\mi{j}^\delta}^2}}\pa{\sum_{\mi{j}\in\mathcal{J}_{\mi{\ell}}}\norm{\mi{j}^\delta}^2\theta_{\mi{j}}(g)^2}}\nonumber\\
	&\leq &\sum_{\mi{\ell}\in \{0,1,\ldots,L\}^d}\sqrt{\sum_{\mi{j}\in\mathcal{J}_{\mi{\ell}}}\cro{\hat\theta_{\mi{j}} - \theta_{\mi{j}}(f)}^2\frac{1}{\norm{\mi{j}^\delta}^2}}, \label{eq:majprop3}
\end{eqnarray}
and this for all $g$ in $\mathcal{W}^\delta(1).$
Since the upper bound does not depend on $g$, we deduce that 
\begin{eqnarray*}
	\esp{d_{\mathcal{W}^{\delta}(1)}\pa{\hat{f}_J,f_J}} &\leq& \sum_{\mi{\ell}\in \{0,1,\ldots,L\}^d}\esp{\sqrt{\sum_{\mi{j}\in\mathcal{J}_{\mi{\ell}}}\cro{\hat\theta_{\mi{j}} - \theta_{\mi{j}}(f)}^2\frac{1}{\norm{\mi{j}^\delta}^2}}} \\
	&\leq& \sum_{\mi{\ell}\in \{0,1,\ldots,L\}^d}\sqrt{\sum_{\mi{j}\in\mathcal{J}_{\mi{\ell}}}\esp{\pa{\hat\theta_{\mi{j}} - \theta_{\mi{j}}(f)}^2}\frac{1}{\norm{\mi{j}^\delta}^2}},
\end{eqnarray*}
by Jensen's inequality. 

Yet, for any $\mi{\ell}$ in $\ac{0,1,\ldots,L}^d$ and any $\mi{j}$ in $\mathcal{J}_{\mi{\ell}}$, by Proposition \ref{prop:CBmechanismespvar}, for all $1\leq i\leq n$, $\esp{Z_{i,\mi{j}}\middle|X_i}=\varphi_{\mi{j}}(X_i)$, thus
$$
\esp{\hat{\theta}_{\mi{j}}} = \frac{1}{n}\sum_{i=1}^n\esp{Z_{i,\mi{j}}}
= \frac{1}{n}\sum_{i=1}^n\esp{\varphi_{\mi{j}}(X_i)} = \theta_{\mi{j}}(f).
$$
Therefore, using once again Proposition \ref{prop:CBmechanismespvar}, with $B_0=\sqrt{2}^d$,
\begin{align*}
	\esp{\pa{\hat\theta_{\mi{j}} - \theta_{\mi{j}}(f)}^2} 
	&  =   \var{\hat\theta_{\mi{j}}} = \frac{1}{n^2}\sum_{i=1}^n \var{Z_{i,\mi{j}}} 
	\leq (2\xi_A B_0)^2\frac{d_{\mi{\ell}}}{n\alpha_{\mi{\ell}}^2}.
\end{align*}
This upper bound no longer depends on $\mi{j}$. 

Moreover, recall the inequality of arithmetic and geometric means that is, for all positive numbers $a_1,\ldots,a_d$, 
\begin{equation}\label{eq:comparemoy}
	\frac{1}{d}\sum_{m=1}^d a_m \geq \pa{\prod_{m=1}^d a_m}^{1/d}.
\end{equation}
It is a direct application of Jensen's inequality for the logarithmic function which is concave. \\

Then, we obtain that for all $\mi{j}$ in $\mathcal{J}_{\mi{\ell}}$, 
\begin{equation}\label{eq:minorjdelta}
	\norm{\mi{j}^\delta}^2 = \sum_{m=1}^d j_m^{2\delta} 
	\ \geq\ 
	d\prod_{m=1}^d j_m^{2\delta/d}
	\ \geq \ d \prod_{m=1}^d 2^{\ell_m 2\delta/d}
	\ =\ d \times d_{\mi{\ell}}^{2\delta/d},
\end{equation}
where the last inequality comes from the definition of $\mathcal{J}_{\mi{\ell}}$. 
Therefore, as the cardinal of $\mathcal{J}_{\mi{\ell}}$ equals $d_{\mi{\ell}}$,
$$
\sum_{\mi{j}\in\mathcal{J}_{\mi{\ell}}}\esp{\pa{\hat\theta_{\mi{j}} - \theta_{\mi{j}}(f)}^2}\frac{1}{\norm{\mi{j}^\delta}^2} 
\leq \frac{(2\xi_A B_0)^2}{d n\alpha_{\mi{\ell}}^2} \sum_{\mi{j}\in \mathcal{J}_{\mi{\ell}}} d_{\mi{\ell}}^{1-2\delta/d}
= \frac{(2\xi_A B_0)^2}{d} \frac{d_{\mi{\ell}}^{2(1-\delta/d)}}{n\alpha_{\mi{\ell}}^2}
= (\tau_{A,d})^2\sigma_{\mi{\ell}}^2, 
$$
where $\tau_{A,d}=2\xi_AB_0/\sqrt{d}=2\sqrt{(2^{d}/d)} A(e^A+1)/(e^A-1)$. 
We deduce that 
\begin{equation*}\label{eq:upperE1}
	\esp{d_{\mathcal{W}^{\delta}(1)}\pa{\hat{f}_J,f_J}} \leq \tau_{A,d}\sum_{\mi{\ell}\in \{0,1,\ldots,L\}^d} \sigma_{\mi{\ell}}, 
\end{equation*} 
which ends the proof of Proposition \ref{prop:controlvar}. 

\subsection{Proof of Lemma \ref{lm:controlbias}}
\label{sect:proof:lm:controlbias}

In this proof, denote for all $\mi{j}=(j_1,\ldots,j_d)$ in $(\N^*)^d$, $\norm{\mi{j}^\beta}^2 = j_1^{2\beta} + j_2^{2\beta} + \ldots +j_d^{2\beta}$, and let $f=\sum_{\mi{j}\in(\N^*)^d}\theta_{\mi{j}}(f) \varphi_{\mi{j}}$ in $\mathcal{W}^{\beta}(R)$, such that 
$$
\sum_{\mi{j}\in(\N^*)^d} \norm{\mi{j}^\beta}^2\theta_{\mi{j}}(f)^2\leq R^2.
$$
We aim at upper bounding 
$$
d_{\mathcal{W}^{\delta}(1)}\pa{f_J,f} = \sup_{g\in\mathcal{W}^\delta(1)}\int_{[0,1]^d}(f_J - f)g.
$$
Let $g$ in $\mathcal{W}^\delta(1)$ and consider $g=\sum_{\mi{j}\in(\N^*)^d}\theta_{\mi{j}}(g) \varphi_{\mi{j}}$ its decomposition in the Fourier basis. 
Then, by the Cauchy-Schwarz inequality, 
$$
\int_{[0,1]^d}(f_J - f)g = \sum_{\mi{j}\notin\{1,\ldots,J\}^d} \cro{-\theta_{\mi{j}}(f)}\theta_{\mi{j}}(g)
\leq \sqrt{\pa{\sum_{\mi{j}\notin\{1,\ldots,J\}^d} \theta_{\mi{j}}(f)^2} \pa{\sum_{\mi{j}\notin\{1,\ldots,J\}^d} \theta_{\mi{j}}(g)^2}}.
$$
Yet, for all $\mi{j}\notin\{1,2,\ldots,J\}^d$, there exists $m$ in $\{1,\ldots,d\}$ such that $j_m>J$. Then, 
$$
\norm{\mi{j}^\beta}^2 \geq j_m^{2\beta} > J^{2\beta}.
$$
Therefore, since $f$ belongs to $\mathcal{W}^\beta(R)$, 
$$
\sum_{\mi{j}\notin\{1,\ldots,J\}^d} \theta_{\mi{j}}(f)^2 = \sum_{\mi{j}\notin\{1,\ldots,J\}^d} \norm{\mi{j}^\beta}^2\theta_{\mi{j}}(f)^2 \frac{1}{\norm{\mi{j}^\beta}^2} \leq R^2 J^{-2\beta}, 
$$
In the same way, since $g$ belongs to $\mathcal{W}^\delta(1)$, 
$\sum_{\mi{j}\notin\{1,\ldots,J\}^d} \theta_{\mi{j}}(g)^2 \leq J^{-2\delta}.$
Hence, 
$$
\int_{[0,1]^d}(f_J - f)g \leq \sqrt{R^2 J^{-2\beta}J^{-2\delta}}.
$$
This being true for any $g$ in $\mathcal{W}^\delta(1)$, we deduce that 
\begin{equation*}\label{eq:upperE2}
	d_{\mathcal{W}^{\delta}(1)}\pa{f_J,f} \leq R J^{-(\beta+\delta)}.
\end{equation*}

\subsection{Proof of Theorem \ref{thm:upper}}
\label{sect:proof:thm:upper}

Applying Proposition \ref{prop:controlvar} and Lemma \ref{lm:controlbias} in the bias-variance decomposition \eqref{eq:biaisvardecomp} leads to 
$$
\esp{d_{\mathcal{W}^{\delta}(1)}\pa{\hat{f}_J,f}} 
\leq \tau_{A,d}\Sigma_J + R J^{-(\beta+\delta)}.
$$
Moreover, for the particular choice of block privacy levels defined in Equation \eqref{eq:defalphaell}, the variance term satisfies 
$$
\Sigma_J = \frac{S_{d,\delta}(J)^2}{\sqrt{n\alpha^2}},
$$
as said in Equation \eqref{eq:defSigmaJ}, where, by definition of $S_{d,\delta}(J)$, and of $d_{\mi{\ell}}$ in \eqref{eq:defdell}, 
$$
S_{d,\delta}(J) 
= \sum_{\mi{\ell}\in\{0,1,\ldots,L\}^d} \prod_{m=1}^d 2^{\ell_m(1-\delta/d)/2}
= \pa{\sum_{\ell=0}^L \cro{2^{(1-\delta/d)/2}}^\ell}^d.
$$
Note that if $\delta\neq d$, then, 
$$
S_{d,\delta}(J)^2 
= \pa{\frac{2^{(1-\delta/d)(L+1)/2}-1}{2^{(1-\delta/d)/2}-1}}^{2d} \\
= \frac{1}{\cro{2^{(1-\delta/d)/2}-1}^{2d}} \pa{(J+1)^{(1-\delta/d)/2}-1}^{2d}.
$$
From here, we distinguish three cases. 
\begin{description}
	\item[Case 1:]
	If $\delta<d$, then $S_{d,\delta}(J)^2 \leq C_{d,\delta} J^{d-\delta}$ and 
	$
	\Sigma_J \leq \frac{C_{d,\delta}}{\sqrt{n\alpha^2}} J^{d-\delta}.
	$
	Thus, the variance and the bias are of the same order if 
	$
	J^{-(\beta+\delta)} \asymp \frac{1}{\sqrt{n\alpha^2}} J^{d-\delta}
	$, that is if $
	J \asymp \pa{n\alpha^2}^{\frac{1}{2\beta+2d}}.
	$
	In that case, we obtain
	$$
	\esp{d_{\mathcal{W}^{\delta}(1)}\pa{\hat{f}_J,f}} \leq C_{A,d,\delta,\beta,R} \pa{n\alpha^2}^{\frac{-(\beta+\delta)}{2\beta+2d}}.
	$$
	\item[Case 2:] 
	If $\delta>d$, then $S_{d,\delta}(J)^2 \leq C_{d,\delta}$ and 
	$
	\Sigma_J \leq \frac{C_{d,\delta}}{\sqrt{n\alpha^2}}.
	$
	Thus, the variance and the bias are of the same order if 
	$
	J^{-(\beta+\delta)} \asymp \frac{1}{\sqrt{n\alpha^2}}
	,$ that is if $
	J \asymp \pa{n\alpha^2}^{\frac{1}{2\beta+2\delta}}.
	$
	In that case,  we obtain
	$$
	\esp{d_{\mathcal{W}^{\delta}(1)}\pa{\hat{f}_J,f}} \leq C_{A,d,\delta,\beta,R} \pa{n\alpha^2}^{-1/2}.
	$$
	
	\item[Case 3:] If $\delta = d$, then $S_{d,\delta}(J)^2 = (L+1)^{2d} \leq C_{d} \log(J)^{2d}$, as $J=2^{L+1}-1$. 
	Then, setting 
	$$
	J \asymp \pa{\frac{n\alpha^2}{\log(n\alpha^2)^{4d}}}^{\frac{1}{2\beta+2d}},
	$$
	leads to 
	$$
	\Sigma_J \leq \frac{C_{d}}{\sqrt{n\alpha^2}}\log(J)^{2d}
	\leq \frac{C_{d,\beta}}{\sqrt{n\alpha^2}}\log(n\alpha^2)^{2d},
	$$
	since 
	$
	\log(J) = \frac{1}{2\beta+2d}\cro{\log(n\alpha^2) - 4d\log(\log(n\alpha^2))} \leq C_{d,\beta}\log(n\alpha^2).
	$
	In that case, we obtain 
	$$
	\esp{d_{\mathcal{W}^{\delta}(1)}\pa{\hat{f}_J,f}} \leq C_{A,d,\beta,R} \pa{\frac{n\alpha^2}{\log(n\alpha^2)^{4d}}}^{-1/2}.
	$$
	
\end{description}

\subsection{Proof of Lemma \ref{lm:ineq}}
\label{sect:proof:lm:ineq}

This proof is based on Hoeffding's inequality (see \cite{hoeffding1963probability}). Fix $J = 2^{L+1}-1$ in $\mathcal{M}_n^\alpha$. 
By taking the supremum for $g$ in $\mathcal{W}^\delta(1)$ in Equation \eqref{eq:majprop3} of the proof of Proposition \ref{prop:controlvar}, we deduce that 
\begin{eqnarray*}
	d_{\mathcal{W}^{\delta}(1)}\pa{\hat f_{J}, f_{J}} 
	&\leq &  \sum_{\mi{\ell}\in \{0,1,\ldots,L\}^d}\sqrt{\sum_{\mi{j}\in\mathcal{J}_{\mi{\ell}}}\cro{\hat\theta_{\mi{j}} - \theta_{\mi{j}}(f)}^2\frac{1}{\norm{\mi{j}^\delta}^2}}.
\end{eqnarray*}
Fix $\mi{\ell}$ in $\{0,\ldots,L\}^d$. For all $\mi{j}$ in $\mathcal{J}_{\mi{\ell}}$, note that, by Proposition \ref{prop:CBmechanismespvar}, 
$$\hat\theta_{\mi{j}} - \theta_{\mi{j}}(f)=  \frac{1}{n}\sum_{i=1}^{n} \left( Z_{i, \mi{j}}-\esp{Z_{i, \mi{j}}}\right),$$  
where the random variables $(Z_{i, \mi{j}})_{1\leq i\leq n}$ are i.i.d. and satisfy
$\abs{Z_{i,\mi{j}}} \leq (2\xi_A B_0) \frac{\sqrt{d_{\mi{\ell}}}}{\alpha_{\mi{\ell}}},$ with $B_0=2^{d/2}$. 
We may thus apply Hoeffding's inequality, and we obtain for all positive $t$, 
\begin{eqnarray*}
	\proba{\left \vert\hat\theta_{\mi{j}} - \theta_{\mi{j}}(f) \right\vert \geq t}\leq 2\exp\left(-\frac{n \alpha_{\mi{\ell}}^2}{ 2(2\xi_A B_0)^2 d_{\mi{\ell}}}t^2 \right).
\end{eqnarray*}

\noindent By setting $\displaystyle t^2=2(2\xi_A B_0)^2\frac{d_{\mi{\ell}}}{n\alpha_{\mi{\ell}}^2}(x_J+x)$, for $x>0$ and $x_J>0$ calibrated below, it follows that
\begin{eqnarray*}
	\proba{\pa{\hat\theta_{\mi{j}} - \theta_{\mi{j}}(f)}^2 \geq 2(2\xi_A B_0)^2\frac{d_{\mi{\ell}}}{n\alpha_{\mi{\ell}}^2}(x_J+x)}
	\leq 2\exp\pa{-(x_J+x)}. 
\end{eqnarray*}

\noindent Moreover, recall that by Equation \eqref{eq:minorjdelta}, $\norm{\mi{j}^\delta}^2\geq d d_{\mi{\ell}}^{2\delta/d}$. Hence, since by Equation \eqref{eq:defsigmaell}, $\tau_{A,d}=2\sqrt{(2^{d}/d)} A(e^A+1)/(e^A-1)=2\xi_AB_0/\sqrt{d}$,
\begin{eqnarray*}
	\proba{\frac{\pa{\hat\theta_{\mi{j}} - \theta_{\mi{j}}(f) }^2}{\norm{\mi{j}^\delta}^2} 
		\geq 2\tau_{A,d}^2\frac{d_{\mi{\ell}}^{1-2\delta/d}}{n\alpha_{\mi{\ell}}^2}(x_J+x)}
	\leq 2\exp\left(-(x_J+x)\right).
\end{eqnarray*}
Thus, by summing over $\mi{j}$ in $\mathcal{J}_{\mi{\ell}}$ (with cardinal $d_{\mi{\ell}}$) and applying the square root, we obtain
$$
\proba{\sqrt{\sum_{\mi{j}\in\mathcal{J}_{\mi{\ell}}}\pa{\hat\theta_{\mi{j}} - \theta_{\mi{j}}(f)}^2\frac{1}{\norm{\mi{j}^\delta}^2}}
	\geq \sqrt{2}\tau_{A,d}  \frac{d_{\mi{\ell}}^{1-\delta/d}}{\sqrt{n}\alpha_{\mi{\ell}}} \sqrt{x_J+x}
}
\leq  2 d_{\mi{\ell}} e^{-x_J} e^{-x}.
$$
Hence, using the definition of $\sigma_{\mi{\ell}}$ and $\Sigma_{J}$ in Equation \eqref{eq:defsigmaell}, 
summing over $\mi{\ell}$ in $\{0,\ldots,L\}^d$ leads to
$$
\proba{
	d_{\mathcal{W}^{\delta}(1)}\pa{\hat f_{J}, f_{J}} 
	\geq \sqrt{2}\tau_{A,d}\Sigma_J(\sqrt{x_J}+\sqrt{x})
}
\leq  2 \pa{\sum_{\mi{\ell}\in \{0,1,\ldots,L\}^d} d_{\mi{\ell}} }e^{-x_J} e^{-x}
= 2\pa{J^d e^{-x_J}} e^{-x}, 
$$
since for all positive $u$ and $v$, $\sqrt{u+v}\leq \sqrt{u} + \sqrt{v}$.

\noindent Now let calibrate $x_J$ such that the probability is controlled. 
Let 
$$x_J = d\log(J) + \frac{3}{2}\log\pa{n\alpha^2} + \log\pa{\tau_{A,d}\Sigma_J}
\quad \text{such that}\quad 
J^d e^{-x_J} = \frac{1}{(n\alpha^2)^{3/2}\tau_{A,d}\Sigma_J}.$$
Finally, recall that 
$V(J) = \sqrt{2}\tau_{A,d}\Sigma_J\sqrt{x_J}.$
Hence, for $\displaystyle t=\sqrt{2}\tau_{A,d}\Sigma_J\sqrt{x},$ we obtain 
$$\proba{d_{\mathcal{W}^{\delta}(1)}\pa{\hat{f}_{J},f_J} \geq V(J) + t} 
\leq \frac{1}{(n\alpha^2)^{3/2}}\frac{2}{\tau_{A,d}\Sigma_J}\exp\pa{\frac{-t^2}{2(\tau_{A,d}\Sigma_J)^2}},$$
which ends the proof. 

\subsection{Proof of Theorem \ref{thm:oraclebound}}
\label{sect:proof:thm:oraclebound}

Fix $J_0 = 2^{L_0+1}-1$ in $\mathcal{M}_n^\alpha$. Then, 
$$d_{\mathcal{W}^{\delta}(1)}\pa{\hat{f}_{\hat{J}},f} 
\leq  d_{\mathcal{W}^{\delta}(1)}\pa{\hat{f}_{\hat{J}},\hat{f}_{\hat{J}\wedge J_0}} + d_{\mathcal{W}^{\delta}(1)}\pa{\hat{f}_{\hat{J}\wedge J_0},\hat{f}_{J_0}} + d_{\mathcal{W}^{\delta}(1)}\pa{\hat{f}_{J_0},f}. 
$$
Yet, 
\begin{eqnarray*}
	d_{\mathcal{W}^{\delta}(1)}\pa{\hat{f}_{\hat{J}},\hat{f}_{\hat{J}\wedge J_0}}
	&\leq & \pa{d_{\mathcal{W}^{\delta}(1)}\pa{\hat{f}_{\hat{J}},\hat{f}_{\hat{J}\wedge J_0}} - 2V(\hat{J})}_+ + 2V(\hat{J}) \\
	&\leq & \hat{A}(J_0) + 2V(\hat{J}).
\end{eqnarray*}
Similarly, 
$$
d_{\mathcal{W}^{\delta}(1)}\pa{\hat{f}_{\hat{J}\wedge J_0},\hat{f}_{J_0}} \leq \hat{A}(\hat{J}) + 2V(J_0).
$$
Therefore, as $\hat{J}$ minimizes the criterion $\widehat{\operatorname{Crit}}(J)$, 
\begin{eqnarray*}
	d_{\mathcal{W}^{\delta}(1)}\pa{\hat{f}_{\hat{J}},f} 
	&\leq & \hat{A}(J_0) + 2V(\hat{J}) + \hat{A}(\hat{J}) + 2V(J_0) + d_{\mathcal{W}^{\delta}(1)}\pa{\hat{f}_{J_0},f} \\
	&\leq & 2\hat{A}(J_0) + 4 V(J_0) + d_{\mathcal{W}^{\delta}(1)}\pa{\hat{f}_{J_0},f}, 
\end{eqnarray*}
and we obtain that 
\begin{equation}\label{eq:1ercontroloracle}
	\esp{d_{\mathcal{W}^{\delta}(1)}\pa{\hat{f}_{\hat{J}},f} }
	\leq 2\esp{\hat{A}(J_0)} + 4 V(J_0) + \esp{d_{\mathcal{W}^{\delta}(1)}\pa{\hat{f}_{J_0},f}}. 
\end{equation}
Let us now control $\esp{\hat{A}(J_0)}$. \\
Note that for all $J\leq J_0$, $\pa{d_{\mathcal{W}^{\delta}(1)}  
	\pa{\hat f_{J},\hat f_{J\wedge J_0}}-2V(J)}_{+}=0$. Therefore, 
$$
\hat{A} (J_0) 
=\max_{\underset{J\geq J_0}{J\in \mathcal{M}_n^\alpha}}\ac{\pa{d_{\mathcal{W}^{\delta}(1)} \pa{\hat f_{J},\hat f_{J_0}}-2V(J)}_{+}}. 
$$
Yet, for all $J\geq J_0$ in $\mathcal{M}_n^\alpha$, $V(J) \geq V(J_0)$. Hence, 
\begin{eqnarray*}
	d_{\mathcal{W}^{\delta}(1)} \pa{\hat f_{J},\hat f_{J_0}}-2V(J)
	&\leq & \pa{d_{\mathcal{W}^{\delta}(1)} \pa{\hat f_{J},f_{J}}-V(J)}_{+} \\ 
	&& +\ d_{\mathcal{W}^{\delta}(1)}\pa{f_{J},f_{J_0}} \\
	&& +\ \pa{d_{\mathcal{W}^{\delta}(1)} \pa{\hat f_{J_0},f_{J_0}}-V(J_0)}_{+}. 
\end{eqnarray*}
Moreover, since 
$$d_{\mathcal{W}^{\delta}(1)}\pa{f_{J},f_{J_0}} \leq 2 \max_{\underset{J\geq J_0}{J \in \mathcal{M}_n^\alpha}} d_{\mathcal{W}^{\delta}(1)}\pa{f_{J},f},$$
we obtain that 
$$\hat{A} (J_0)
\leq 2 \max_{\underset{J\geq J_0}{J \in \mathcal{M}_n^\alpha}} \ac{\pa{d_{\mathcal{W}^{\delta}(1)} \pa{\hat f_{J},f_{J}}-V(J)}_{+} } 
+ 2 \max_{\underset{J\geq J_0}{J \in \mathcal{M}_n^\alpha}} d_{\mathcal{W}^{\delta}(1)}\pa{f_{J},f}.$$
Denote for simplicity $D_J = d_{\mathcal{W}^{\delta}(1)} \pa{\hat f_{J},f_{J}}-V(J).$
Then, 
\begin{equation}
	\label{eq:majespAhat}
	\esp{\hat{A} (J_0)} \leq 2 \esp{\max_{\underset{J\geq J_0}{J \in \mathcal{M}_n^\alpha}} \pa{D_J}_{+} } + 2 \max_{\underset{J\geq J_0}{J \in \mathcal{M}_n^\alpha}} d_{\mathcal{W}^{\delta}(1)}\pa{f_{J},f}.
\end{equation}
We use the concentration inequality of Lemma \ref{lm:ineq} to control the first term in the upper-bound above. 
First notice that 
$$\esp{\max_{\underset{J\geq J_0}{J \in \mathcal{M}_n^\alpha}} \pa{D_J}_{+} }
=\int_{0}^{+\infty} \proba{\max_{\underset{J\geq J_0}{J \in \mathcal{M}_n^\alpha}} \pa{D_J}_{+}>t}dt. 
$$
Moreover, for all positive $t$, 
\begin{eqnarray*}
	\proba{\max_{\underset{J\geq J_0}{J \in \mathcal{M}_n^\alpha}} \pa{D_J}_{+} >t}
	&\leq & \sum_{\underset{J\geq J_0}{J \in \mathcal{M}_n^\alpha}} \proba{\pa{D_J}_{+}>t}.
\end{eqnarray*}
Hence, 
$$\esp{\max_{\underset{J\geq J_0}{J \in \mathcal{M}_n^\alpha}} \ac{\pa{D_J}_{+} }}
\leq \sum_{\underset{J\geq J_0}{J \in \mathcal{M}_n^\alpha}} \int_{0}^{+\infty} \proba{\pa{D_J}_{+}>t} dt. $$
Yet, by Lemma \ref{lm:ineq}, for all positive $t$,
\begin{eqnarray*}
	\proba{\pa{D_J}_{+}>t} \ \leq\ \proba{D_J>t} 
	&=& \proba{d_{\mathcal{W}^{\delta}(1)} \pa{\hat{f}_{J},f_{J}} > V(J) + t}\\
	&\leq & \frac{\sqrt{2\pi}}{(n\alpha^2)^{3/2}}\times\frac{2}{\sqrt{2\pi}(\tau_{A,d}\Sigma_J)}\exp\pa{\frac{-t^2}{2(\tau_{A,d}\Sigma_J)^2}}.
\end{eqnarray*}
Therefore, by integrating on $\R_+$, and summing over $J\geq J_0$ in $\mathcal{M}_n^\alpha$, we obtain 
$$\esp{\max_{\underset{J\geq J_0}{J \in \mathcal{M}_n^\alpha}} \ac{\pa{D_J}_{+} }}
\leq \sqrt{2\pi}\frac{\abs{\mathcal{M}_n^\alpha}}{\sqrt{n\alpha^2}}\frac{1}{n\alpha^2}.$$
Furthermore, since $2^{L+1}-1\leq n\alpha^2$ if and only if $L \leq \log_2(n\alpha^2+1)-1$, then
$$\abs{\mathcal{M}_n^\alpha} = \lfloor\log_2(n\alpha^2+1)\rfloor \leq C\log(n\alpha^2),$$
where $C$ is a universal constant. Moreover, $\log(n\alpha^2)\leq 2e^{-1}\sqrt{n\alpha^2}$. Hence, we deduce that $\abs{\mathcal{M}_n^\alpha}\leq C\sqrt{n\alpha^2}$, and therefore, by \eqref{eq:majespAhat},
\begin{equation}\label{eq:controlespAhatordrebiais}
	\esp{\hat{A} (J_0)} \leq 2 \max_{\underset{J\geq J_0}{J \in \mathcal{M}_n^\alpha}} d_{\mathcal{W}^{\delta}(1)}\pa{f_{J},f} + \frac{C}{n\alpha^2}. 
\end{equation}
Finally, using \eqref{eq:1ercontroloracle}, we obtain that 
$$\esp{d_{\mathcal{W}^{\delta}(1)}\pa{\hat{f}_{\hat{J}},f}}
\leq \ac{\esp{d_{\mathcal{W}^{\delta}(1)}\pa{\hat{f}_{J_0},f}} + 4 \max_{\underset{J\geq J_0}{J \in \mathcal{M}_n^\alpha}} d_{\mathcal{W}^{\delta}(1)}\pa{f_{J},f} + 4 V(J_0)} + \frac{C}{n\alpha^2}, 
$$
and this for all $J_0$ in $\mathcal{M}_n^\alpha$.

\subsection{Proof of Corollary \ref{cor:adaptupper}}
\label{sect:proof:cor:adaptupper}

Recall that, from Equation \eqref{eq:biaisvaradapt}, the adaptive estimator $\hat{f}_{\hat{J}}$ satisfies 
$$
\esp{d_{\mathcal{W}^{\delta}(1)}\pa{\hat{f}_{ \hat{J}},f}}
\leq C\inf_{J_0\in \mathcal{M}_n^\alpha} \pa{\max_{\underset{J\geq J_0}{J\in \mathcal{M}_n^\alpha}}d_{\mathcal{W}^{\delta}(1)}  \pa{f_{J},f}+ V(J_0)}+\frac{C}{n\alpha^2}. 
$$
As the term proportional to $(n\alpha^2)^{-1}$ is negligible compared to the first one, the aim is to determine a $J_0$ in $\mathcal{M}_n^\alpha$ that allows the compromise between the bias term $\max_{\underset{J\geq J_0}{J\in \mathcal{M}_n^\alpha}}d_{\mathcal{W}^{\delta}(1)}  \pa{f_{J},f}$ and the variance term $V(J_0)$. 
Moreover, as $f$ belongs to $\mathcal{W}^{\beta}(R)$, from Lemma \ref{lm:controlbias}, for all $J\geq J_0$, 
$d_{\mathcal{W}^{\delta}(1)}\pa{f_J,f} \leq R J^{-(\beta+\delta)}.$ We deduce that 
$$\max_{\underset{J\geq J_0}{J\in \mathcal{M}_n^\alpha}}d_{\mathcal{W}^{\delta}(1)}  \pa{f_{J},f} \leq R J_0^{-(\beta+\delta)}.$$
As in the proof of Theorem \ref{thm:upper}, we distinguish three cases to control the variance term 
$$V(J_0) = \sqrt{2}\tau_{A,d}\Sigma_{J_0}\sqrt{d\log(J_0) + \frac{3}{2}\log(n\alpha^2) + \log(\tau_{A,d}\Sigma_{J_0})}.$$ 
The control of $\Sigma_{J_0}$ comes from the choice of the block privacy levels in Equation \eqref{eq:defalphaell} and is exactly the same as in the proof of Theorem \ref{thm:upper} (see Section \ref{sect:proof:thm:upper}). 
\begin{description}
	\item[Case 1:]
	If $\delta<d$, then for all $J_0$ in $\mathcal{M}_n^\alpha$, 
	$
	\Sigma_{J_0} \leq \frac{C_{d,\delta}}{\sqrt{n\alpha^2}} J_0^{d-\delta}.
	$
	In particular, we obtain that $\log(\Sigma_{J_0}) \leq C_{d,\delta}\log(J_0) \leq C_{d,\delta}\log(n\alpha^2)$, and 
	$$V(J_0) \leq C_{A,d,\delta}J_0^{d-\delta}\sqrt{\frac{\log(n\alpha^2)}{n\alpha^2}}.$$
	Thus, $V(J_0)$ and the bias $RJ_0^{-(\beta+\delta)}$ are of the same order if 
	$
	J_0^{-(\beta+\delta)} \asymp \sqrt{\frac{\log(n\alpha^2)}{n\alpha^2}} J_0^{d-\delta}
	$, that is if $
	J_0 \asymp \pa{\frac{n\alpha^2}{\log(n\alpha^2)}}^{\frac{1}{2\beta+2d}}.
	$
	In that case, we obtain
	$$
	\esp{d_{\mathcal{W}^{\delta}(1)}\pa{\hat{f}_{\hat{J}},f}} 
	\leq C_{A,d,\delta,\beta,R} \pa{\frac{n\alpha^2}{\log(n\alpha^2)}}^{\frac{-(\beta+\delta)}{2\beta+2d}}.
	$$
	\item[Case 2:] 
	If $\delta>d$, then 
	$
	\Sigma_{J_0} \leq \frac{C_{d,\delta}}{\sqrt{n\alpha^2}}.
	$
	Thus, as $n\alpha^2\geq 2$, $\log(\Sigma_{J_0}) \leq C_{d,\delta} \leq C_{d,\delta}\log(n\alpha^2)$, and
	$$V(J_0) \leq C_{A,d,\delta}\sqrt{\frac{\log(n\alpha^2)}{n\alpha^2}}.$$
	Hence, $V(J_0)$ and the bias $RJ_0^{-(\beta+\delta)}$ are of the same order if 
	$J_0^{-(\beta+\delta)} \asymp \sqrt{\frac{\log(n\alpha^2)}{n\alpha^2}},$ that is if 
	$J_0 \asymp \pa{\frac{n\alpha^2}{\log(n\alpha^2)}}^{\frac{1}{2\beta+2\delta}}.$
	In that case,  we obtain
	$$
	\esp{d_{\mathcal{W}^{\delta}(1)}\pa{\hat{f}_{\hat{J}},f}} \leq C_{A,d,\delta,\beta,R} \pa{\frac{n\alpha^2}{\log(n\alpha^2)}}^{-1/2}.
	$$
	
	\item[Case 3:] If $\delta = d$, then 
	$
	\Sigma_{J_0} \leq \frac{C_{d}}{\sqrt{n\alpha^2}}\log(J_0)^{2d}
	\leq \frac{C_{d}}{\sqrt{n\alpha^2}}\log(n\alpha^2)^{2d}. 
	$
	In particular, we obtain that $\log(\Sigma_{J_0})\leq C_{d}\log(n\alpha^2)$ and 
	$$
	V(J_0) \leq C_{A,d}\sqrt{\frac{\pa{\log(n\alpha^2)}^{4d+1}}{n\alpha^2}}.
	$$
	Hence, $V(J_0)$ and the bias $RJ_0^{-(\beta+\delta)}$ are of the same order if 
	$J_0^{-(\beta+\delta)} \asymp \sqrt{\frac{\pa{\log(n\alpha^2)}^{4d+1}}{n\alpha^2}},$ that is if 
	$J_0 \asymp \pa{\frac{n\alpha^2}{\pa{\log(n\alpha^2)}^{4d+1}}}^{\frac{1}{2\beta+2\delta}}.$
	In that case,  we obtain
	$$
	\esp{d_{\mathcal{W}^{\delta}(1)}\pa{\hat{f}_{\hat{J}},f}} \leq C_{A,d,\beta,R} \pa{\frac{n\alpha^2}{\pa{\log(n\alpha^2)}^{4d+1}}}^{-1/2}.
	$$
\end{description}
Note that in all cases, the $J_0$ that satisfies the trade-off belongs to the collection $\mathcal{M}_n^\alpha$.



\subsection*{Acknowledgement}
Our work has benefitted from the AI Interdisciplinary Institute ANITI. ANITI is funded by the France 2030 program under the Grant agreement n\textsuperscript{o} ANR-23-IACL-0002. We also recognize the funding by ANITI ANR-19-PI3A-0004, and by CIMI ANR-11-LABX-0040.

\newpage
\appendix

\section{Complementary proofs}\label{sect:complproofs}
\subsection{Characterization of Sobolev balls using partial derivatives}
\label{sect:sobolevcharact}

\begin{lemma}\label{lm:derivtocoefSobolev}
	Let $\mi{\beta}=(\beta_1,\ldots,\beta_d)$ be a regularity parameter with integer coordinates. 
	Consider $f$ a function in  $\L_2([0,1]^d)$ such that  $ \int_{[0,1]^d} f^2(x) dx \leq C_1^2$ and for all $1\leq m\leq d$, $f$ is $\beta_m$-differentiable w.r.t. the $m$-th variable, and satisfies, 
	$$\int_{[0,1]^d}\cro{\frac{\partial^{\beta_m}f}{\partial x_m^{\beta_m}} (x)}^2 dx \leq C_2^2.$$
	Assume also that $f$ and all its partial derivatives  are periodic, namely for all $m =1, \ldots, d$, for all $ 0 \leq l_m \leq \beta _m $, for all $(x_1,\ldots,x_{m-1},x_{m+1},\ldots,x_d)$ in $[0,1]^{d-1}$, 
	$$ \frac{\partial^{l_m}f}{\partial x_m^{l_m}} (x_1, \ldots, x_{m-1},0, \ldots, x_d)= \frac{\partial^{l_m}f}{\partial x_m^{l_m}} (x_1, \ldots, x_{m-1},1, \ldots, x_d).$$
	Then, 
	$$\sum_{(j_1,\ldots,j_d)\in(\N^*)^d} \pa{j_1^{2\beta_1} + \ldots + j_d^{2\beta_d}}\theta_{(j_1,\ldots,j_d)}^2(f)\leq d\pa{C_1^2 + C_2^2}.$$
	In particular, $f$ belongs to the Sobolev ball $\mathcal{W}^{\mi{\beta}}\pa{R}$ as defined in \eqref{eq:defSobolev}, as soon as 
	$$d \pa{C_1^2+C_2^2}\leq R^2.$$
\end{lemma}
First note that, by Parseval's equality,
\begin{equation}\label{eq:Parseval1}
	\sum_{\mi{j} \in (\N^*)^d} \theta^2_{\mi{j}}\pa{f} = \int_{[0,1]^d} f^2(x) dx \leq C_1^2 ,
\end{equation}
and
\begin{equation}\label{eq:Parseval}
	\sum_{\mi{j} \in (\N^*)^d} \theta^2_{\mi{j}}\pa{ \frac{\partial^{\beta_m}f}{\partial x_m^{\beta_m}} } = \int_{[0,1]^d}\cro{\frac{\partial^{\beta_m} f}{\partial x_m^{\beta_m}} (x)}^2 dx \leq C_2^2.
\end{equation}
By integrations by parts, and using the periodicity of $f $ and all its partial derivatives, we have 
\begin{eqnarray*}
	\theta_{\mi{j}}\pa{ \frac{\partial^{\beta_m }f}{\partial x_m^{\beta_m}}}  
	&=&  \int_{[0,1]^d}  {\frac{\partial^{\beta_m} f}{\partial x_m^{\beta_m}}} \times \varphi_{\mi{j}}  \\
	&=& \pm  \int_{[0,1]^d} f \times {\frac{\partial^{\beta_m}  \varphi_{\mi{j}}}{\partial x_m^{\beta_m}}} .
\end{eqnarray*}
Recalling that  for all $\mi{j}=(j_1,\ldots,j_d)$ in $(\N^*)^d$ and for all $x = (x_1,\ldots,x_d)$ in $[0,1]^d$, 
$
\varphi_{\mi{j}}(x) = \prod_{l=1}^{d} \varphi_{j_l}(x_l),
$
we have 
$$ \frac{\partial^{\beta_m} \varphi_{\mi{j}}}{\partial x_m^{\beta_m}}(x)  =  \pa{\prod_{\underset{l\neq m}{l =1}}^{d} \varphi_{j_l}(x_l)} \varphi_{j_m}^{(\beta_m)}(x_m).$$
Moreover, by definition of the Fourier basis, 
\begin{itemize}
	\item if $\beta $ is even, for all $ j \in \N^*$, $ \varphi_{2j}^{(\beta)} = \pm (2 \pi  j )^{\beta}  \varphi_{2j}$ and $ \varphi_{2j+1}^{(\beta)} = \pm (2 \pi  j )^{\beta}  \varphi_{2j+1}$, 
	\item if $\beta $ is odd,  for all $ j \in \N^*$, $ \varphi_{2j}^{(\beta)} = \pm (2 \pi  j )^{\beta}  \varphi_{2j+1}$ and  $ \varphi_{2j+1}^{(\beta)} = \pm (2 \pi  j )^{\beta}  \varphi_{2j}$. 
\end{itemize}
Hence, if $\beta_m$ is even and $j_m \geq 2$, 
$$
\theta_{\mi{j}}\pa{ \frac{\partial^{\beta_m }f}{\partial x_m^{\beta_m}}} =
\left\{\begin{array}{ll}
	\pm \pa{\frac{2 \pi  j_m}2}^{\beta_m} \theta_{\mi{j}} (f) & \mbox{ if } j_m  \mbox{ is even } \\
	\pm \pa{ \frac{2 \pi  (j_m-1)}2}^{\beta_m} \theta_{\mi{j}} (f) & \mbox{ if } j_m  \mbox{ is odd. }
\end{array}\right.
$$
Since for all $j_m  \geq  2 $, $j_m \leq 2(j_m-1)$, we  deduce from \eqref{eq:Parseval} that
\begin{equation}\label{eq:caspair}
	\sum_{\mi{j} \in (\N^*)^d, \ j_m \geq 2}   j_m^{2\beta_m}   \theta^2_{\mi{j}} (f) \leq \sum_{\mi{j} \in (\N^*)^d, \ j_m \geq 2}   \pa{\frac{ \pi  j_m}2}^{2\beta_m}   \theta^2_{\mi{j}} (f) \leq C_2^2.
\end{equation}
If $\beta_m$ is odd and $j_m \geq 2$, 
$$
\theta_{\mi{j}}\pa{ \frac{\partial^{\beta_m }f}{\partial x_m^{\beta_m}}} 
= \left\{\begin{array}{ll}
	\pm \pa{\frac{2 \pi  j_m}2}^{\beta_m} \theta_{j_1\ldots j_{m+1}\ldots j_d} (f) & \mbox{ if } j_m  \mbox{ is even } \\
	\pm \pa{ \frac{2 \pi  (j_m-1)}2}^{\beta_m} \theta_{j_1\ldots j_{m-1}\ldots j_d} (f) &\mbox{ if } j_m  \mbox{ is odd,}
\end{array}\right.
$$
and we also obtain \eqref{eq:caspair}. 
Finally, using \eqref{eq:Parseval1}, we  get
\begin{equation*}
	\sum_{\mi{j} \in (\N^*)^d}   {j_m}^{2\beta_m}   \theta^2_{\mi{j}} (f) \leq C_1^2 + C_2^2,
\end{equation*}
which concludes the proof of Lemma  \ref{lm:derivtocoefSobolev}.

\subsection{Proof of Lemma \ref{lm:inclusionParamdansSobolev}} \label{sect:proofinclusionParamdansSobolev}

\begin{enumerate}
	\item Let $f_\nu$ belong to $\mathcal{F}^\beta(\gamma_n)$, and let $x$ in $[0,1]^d$.
	Then, there exists a unique $\mi{j_0}$ such that $x$ belongs to the support of $G_{\mi{j_0}}$. 
	In particular, if $\gamma_n \leq \norm{\psi}_\infty^{-d} \leq J^{\beta}\norm{\psi}_\infty^{-d},$ then
	$$f_\nu(x) = 1 + \frac{\gamma_n}{J^\beta} \nu_{\mi{j_0}} G_{\mi{j_0}}(x) \geq 1 - \frac{\gamma_n}{J^{\beta}}\abs{G_{\mi{j_0}}(x)} \geq 1 - \frac{\gamma_n}{J^\beta} \norm{\psi}_\infty ^d \geq 0. $$
	
	Moreover, using the left hand side equation in \eqref{eq:momentsGj} directly leads to
	\begin{eqnarray*}
		\int_{[0,1]^d}f_\nu(x)dx 
		&=& 1+\frac{\gamma_n}{J^{\beta}}\sum_{\mi{j}\in\{1,\ldots,J\}^d}\nu_{\mi{j}}\underbrace{\int_{[0,1]^d}G_{\mi{j}}(x)dx}_{=0} 
		= 1.
	\end{eqnarray*}
	
	\item To prove this point, we shall use Lemma \ref{lm:derivtocoefSobolev}. 
	
	Consider $\nu=(\nu_{\mi{j}})_{\mi{j}\in\{1,\ldots,J\}^d}$ in $\{0,1\}^{J^d}$ and consider $f_{\nu}:[0,1]^d\to \R$ in $\mathcal{F}^{\beta}(\gamma_n)$, defined for all $x$ in $[0,1]^d$ by 
	$$
	f_\nu(x) = 1 + \frac{\gamma_n}{J^{\beta}} \sum_{\mi{j}\in\ac{1,\ldots,J}^d} \nu_{\mi{j}} G_{\mi{j}}(x).
	$$
	
	Let us first upper bound the $\L_2$ norm of $f_\nu$. 
	As the supports of the $G_{\mi{j}}$ are disjoint,  
	by both Equations in \eqref{eq:momentsGj}
	\begin{eqnarray*}
		\int_{[0,1]^d}f_\nu(x)^2 dx 
		&=& 1 + \frac{\gamma_n^2}{J^{2\beta}} \sum_{\mi{j}\in\ac{1,\ldots,J}^d} \nu_{\mi{j}}^2 \int_{[0,1]^d}G_{\mi{j}}(x)^2dx.\\
		&\leq & 1 + \frac{\gamma_n^2}{J^{2\beta}} \norm{\psi}_2^{2d} \\
		&\leq & 1 + \gamma_n^2 \norm{\psi}_2^{2d}. 
	\end{eqnarray*}
	
	Let us now upper bound the $\L_2$ norm of the partial derivatives of $f_\nu$. 
	Let $m$ in $\{1,\ldots,d\}$. Then for all $x\in[0,1]^d$,
	$$
	\frac{\partial ^{\beta}}{\partial x_m^\beta} f_{\nu}(x) 
	= \frac{\gamma_n}{J^\beta} \sum_{\mi{j}\in\ac{1,\ldots,J}^d} \nu_{\mi{j}} \cro{\frac{\partial ^{\beta}}{\partial x_m^\beta}G_{\mi{j}}(x)}. 
	$$
	Since the supports of the $G_{\mi{j}}$ for $\mi{j}\in\{1,\ldots,J\}^d$ are disjoint, 
	$$
	\cro{\frac{\partial ^{\beta}}{\partial x_m^\beta} f_{\nu}(x) }^2 
	= \frac{\gamma_n^2}{J^{2\beta}} \sum_{\mi{j}\in\ac{1,\ldots,J}^d} \nu_{\mi{j}}^2 \cro{\frac{\partial ^{\beta}}{\partial x_m^\beta}G_{\mi{j}}(x)}^2, 
	$$
	and in particular, 
	$$
	\int_{[0,1]^d}\cro{\frac{\partial ^{\beta}}{\partial x_m^\beta} f_{\nu}(x) }^2 dx
	= \frac{\gamma_n^2}{J^{2\beta}} \sum_{\mi{j}\in\ac{1,\ldots,J}^d} \nu_{\mi{j}}^2 \int_{[0,1]^d}\cro{\frac{\partial^{\beta}}{\partial x_m^\beta}G_{\mi{j}}(x)}^2dx. 
	$$
	
	Yet, by definition of $G_{\mi{j}}$, 
	\begin{eqnarray*}
		\abs{\frac{\partial^{\beta}}{\partial x_m^\beta}G_{\mi{j}}(x)}
		&=& \abs{J^{\beta}\psi^{(\beta)}\pa{J\pa{x_{m} - \frac{j_{m}-1}{J}}}\times\cro{\prod_{m'\neq m}\psi\pa{J\pa{x_{m'} - \frac{j_{m'}-1}{J}}}}}
	\end{eqnarray*}
	
	Therefore, 
	\begin{eqnarray*}
		\int_{[0,1]^d}\cro{\frac{\partial^{\beta}}{\partial x_m^\beta}G_{\mi{j}}(x)}^2dx 
		&=& J^{2\beta} \int_{\frac{j_{m}-1}{J}}^{\frac{j_m}{J}} \cro{\psi^{(\beta)}\pa{J\pa{x_{m} - \frac{j_{m}-1}{J}}}}^2 dx_m \\
		&& \times \prod_{m'\neq m} \int_{\frac{j_{m'}-1}{J}}^{\frac{j_{m'}}{J}}\cro{\psi\pa{J\pa{x_{m'} - \frac{j_{m'}-1}{J}}}}^2 dx_{m'} \\
		&= & \frac{J^{2\beta}\pa{ \norm{\psi^{(\beta)}}_2 \norm{\psi}_2^{(d-1)}}^2}{J^d}
	\end{eqnarray*}

	We deduce that, as $\nu_{\mi{j}}^2\leq 1$ for all $\mi{j}$, 
	$$
	\int_{[0,1]^d}\cro{\frac{\partial ^{\beta}}{\partial x_m^\beta} f_{\nu}(x) }^2 dx
	\leq \frac{\gamma_n^2}{J^{2\beta}} \times J^d \times \frac{J^{2\beta}\pa{\norm{\psi^{(\beta)}}_2\norm{\psi}_{2}^{d-1}}^2}{J^d} 
	= \gamma_n^2 \pa{\norm{\psi^{(\beta)}}_2\norm{\psi}_{2}^{d-1}}^2.
	$$
	This result being true for all $1\leq m\leq d$, we deduce from Lemma \ref{lm:derivtocoefSobolev} that $f_\nu$ belongs to the isotropic Sobolev ball $\mathcal{W}^\beta(R)$ as soon as 
	$$d\times \ac{1 + \gamma_n^2\norm{\psi}_2^{2(d-1)}  \cro{ \norm{\psi}_2^{2} + \norm{\psi^{(\beta)}}_2^2 }} \leq R^2,$$
	i.e. 
	$$\gamma_n^2\leq \frac{R^2 - d}{d\norm{\psi}_2^{2(d-1)}  \cro{ \norm{\psi}_2^{2} + \norm{\psi^{(\beta)}}_2^2 }} .$$

	The proof is similar for the inclusion $\mathcal{D}^{\delta}(\eta)\subset \mathcal{W}^{\delta}(1),$ as soon as 
	$$d\times \eta^2 \norm{\psi}_2^{2(d-1)}  \cro{ \norm{\psi}_2^{2} + \norm{\psi^{(\delta)}}_2^2 }\leq 1,$$
	which ends the proof. 
\end{enumerate}

\subsection{Proof of Lemma \ref{lm:controlBka}}
\label{sect:proof:lm:controlBka}

Recall that, by definition, 
$$
B_{k}(a) = B_0 \frac{e^{a}+1}{e^{a}-1} \Gamma_{k}, \qquad \text{with}\qquad 
\frac{1}{\Gamma_{k}} =  \frac{1}{2^{k-1}}\binom{k-1}{\lfloor\frac{k-1}{2}\rfloor}. 
$$
On the one hand, let us prove that
\begin{equation}
	\label{eq:majfcta}
	\frac{e^a+1}{e^a-1} \leq \frac{e^A+1}{e^A-1}\times \frac{A}{a} =\frac{\xi_A}{a}.
\end{equation}
To do so, we prove that the function $$g:x \mapsto x\times \frac{e^x+1}{e^x-1}$$
is increasing on $\R_+^*$. Indeed, $g$ is a differentiable function, with derivative 
$$g'(x) = \frac{h(x)}{(e^x-1)^2}
\quad \text{where} \quad 
h(x) = e^{2x} - 2xe^x -1,$$
Moreover, $h$ is also differentiable on $\R$ (and thus on $\R_+^*$), with derivative 
$$h'(x) = 2e^x\cro{e^x - (1+x)} \geq 0.$$
Therefore, $h$ is increasing, and as $h(0) = 0$, $h$ is positive on $\R_+^*$, and so is $g'$. \\
On the other hand, let us prove that 
\begin{equation}
	\label{eq:majGammak}
	\Gamma_k\leq 2\sqrt{k}.
\end{equation}
To do so, we use the inequality proved by \cite{robins1955remark} on Stirling's formula, that implies that for all positive integer $q$, 
$$\sqrt{2\pi q} \pa{\frac{q}{e}}^{q} \leq q! \leq \sqrt{2\pi q}\pa{\frac{q}{e}}^{q} e^{\frac{1}{12 q}}.$$ 
Hence, since for all $q\geq 1$, $e^{1/(6q)} \leq e^{1/6}$, 
\begin{equation}
	\label{eq:controlqparmi2q}
	2^{2q} \frac{(q!)^2}{(2q)!} 
	\leq 2^{2q} \frac{2\pi q (q/e)^{2q}}{\sqrt{4\pi q}(2q/e)^{2q}} e^{2/(12q)} 
	\leq e^{1/6}\sqrt{\pi q}. 
\end{equation}
\begin{itemize}
	\item First notice that $\Gamma_1=1\leq c_1\sqrt{1},$ with $c_1=1\leq 2$. 
	\item Second, notice that $\Gamma_2=2\leq c_2\sqrt{2},$ with $c_2=\sqrt{2}\leq 2.$
	\item Third, if $k=2p+1$ is odd, with $p\geq 1$, then from \eqref{eq:controlqparmi2q}, 
	$$
	\Gamma_k = 2^{2p} \frac{(p!)^2}{(2p)!} 
	\leq e^{1/6}\sqrt{\pi p} 
	\leq e^{1/6}\sqrt{\frac{\pi k}{2}} 
	= c_3\sqrt{k}, 
	$$
	where $c_3 = e^{1/6}\sqrt{\pi/2}\approx 1.48\leq 2$. 
	\item Fourth, if $k=2p$ is even, with $p\geq 2$, since $(2p)/(2p-1) \leq 4/3$, then from \eqref{eq:controlqparmi2q}, 
	$$\Gamma_k = 2^{2p-1} \frac{(p-1)! p!}{(2p-1)!} 
	= \frac{2p}{2p-1} \times 2^{2(p-1)} \frac{((p-1)!)^2}{[2(p-1)]!} 
	\leq \frac{4e^{1/6}}{3}\sqrt{\pi (p-1)} = c_4\sqrt{k},$$
	where $\displaystyle c_4 = \frac{2\sqrt{2\pi}e^{1/6}}{3}\approx 1.97\leq 2$. 
\end{itemize}
Combining all case leads to \eqref{eq:majGammak}, which ends the proof of Lemma \ref{lm:controlBka}.

\section{Why blocks are necessary for minimax optimality}
\label{app:globalnotminimax}

Consider the \emph{Coordinate global privacy mechanism} from \cite{duchi2018minimax} and adapted by \cite{butucea2023phase}. It is exactly the same mechanism where all the coefficients $\pa{\varphi_{\mi{j}}(X_i)}_{\mi{j}\in\{1,2,\ldots,J\}^d}$ are privatized at once. 
In particular, as in Proposition \ref{prop:CBmechanismespvar}, the variance of the private views of $X_1,\ldots,X_n$ satisfies 
\begin{equation}\label{eq:varZglobal}
	\var{Z_{i,\mi{j}}}\leq C_{A,d}^2 \frac{J^d}{\alpha^2}.
\end{equation}
However, the control of the variance term in Proposition \ref{prop:controlvar} differs. 
More precisely, 
\begin{equation}\label{eq:controlvarbis}
	\esp{d_{\mathcal{W}^{\delta}(1)}\pa{\hat{f}_J,f_J}} \leq C_{A,d}\tilde{\Sigma}_J,
\end{equation}
where 
\begin{equation*}\label{eq:controlvarbis2}
	\tilde{\Sigma}_J =  \frac{\pa{J\times \tilde{S}_{\delta,d}(J)}^{d/2}}{\sqrt{n\alpha^2}}
	\quad \text{and} \quad 
	\tilde{S}_{\delta,d}(J) = \sum_{j=1}^J \frac{1}{j^{2\delta/d}}.
\end{equation*}
Indeed, one may easily prove with the same arguments as in the proof of Proposition \ref{prop:controlvar} (see Section \ref{sect:proof:prop:controlvar}) that
$$
\esp{d_{\mathcal{W}^{\delta}(1)}\pa{\hat{f}_J,f_J}} \leq \sqrt{\sum_{\mi{j}\in\{1,2,\ldots,J\}^d}\esp{\pa{\hat\theta_{\mi{j}} - \theta_{\mi{j}}(f)}^2}\frac{1}{\norm{\mi{j}^\delta}^2}},
$$
where this time, using Equation \eqref{eq:varZglobal}, 
$$
\esp{\pa{\hat\theta_{\mi{j}} - \theta_{\mi{j}}(f)}^2} 
=   \var{\hat\theta_{\mi{j}}} 
\leq C_{A,d}^2 \frac{J^d}{n\alpha^2}. 
$$
Hence, it remains to control
\begin{equation*}\label{eq:sumpenible}\sum_{\mi{j}\in\{1,2,\ldots,J\}^d}\frac{1}{\norm{\mi{j}^\delta}^2} = \sum_{(j_1,\ldots,j_d)\in\{1,2,\ldots,J\}^d} \frac{1}{j_1^{2\delta} + \dots + j_d^{2\delta}}.
\end{equation*}
Yet, using the inequality of arithmetic and geometric means \eqref{eq:comparemoy}, we obtain 
$$\norm{\mi{j}^\delta}^2 = \sum_{m=1}^d j_m^{2\delta} 
\ \geq\ 
d\prod_{m=1}^d j_m^{2\delta/d},$$
and thus 
$$\sum_{\mi{j}\in\{1,2,\ldots,J\}^d}\frac{1}{\norm{\mi{j}^\delta}^2} 
\leq \frac{1}{d} \sum_{(j_1,\ldots,j_d)\in\{1,2,\ldots,J\}^d} \cro{\prod_{m=1}^d \frac{1}{j_m^{2\delta/d}}} 
= \frac{1}{d}\pa{\sum_{j=1}^J \frac{1}{j^{2\delta/d}}}^{d} 
= \frac{1}{d} {\tilde{S}_{\delta,d}(J)}^d,$$
which ends the proof of Equation \eqref{eq:controlvarbis}. 
Hence the control of the variance depends on the behaviour of $\tilde{S}_{\delta,d}(J)$, and thus on the value of $2\delta/d$. 

\paragraph{}
In the bias-variance decomposition \eqref{eq:biaisvardecomp}, the control of the bias term in Lemma \ref{lm:controlbias} remains unchanged since it comes from the regularity of the density $f$. 
Therefore, as in the proof of Theorem  \ref{thm:upper}, we obtain
$$
\esp{d_{\mathcal{W}^{\delta}(1)}\pa{\hat{f}_J,f}} \leq C_{A,d}\tilde\Sigma_J + R J^{-(\beta+\delta)}.
$$

As with the block mechanism, we distinguish three cases. However, unlike with the block mechanism, the threshold is when $\delta=d/2$ instead of $\delta=d$. 
\begin{description}
	\item[Case 1:]
	If $\delta<d/2$, then a comparison between sum and integral leads to 
	$$\tilde{S}_{\delta,d}(J) = \sum_{j=1}^J \frac{1}{j^{2\delta/d}} \leq \int_{0}^Jx^{-2\delta/d}dx = \frac{J^{1-2\delta/d}}{1-2\delta/d},$$
	and thus
	$$
	\tilde{\Sigma}_J \leq C_{d,\delta}\frac{J^{d-\delta}}{\sqrt{n\alpha^2}}.
	$$
	The variance and bias terms are of the same order if $J^{-(\beta+\delta)} \asymp \frac{1}{\sqrt{n\alpha^2}} J^{d-\delta}$, that is if 
	$J \asymp \pa{n\alpha^2}^{\frac{1}{2\beta+2d}}.$
	In that case, we obtain
	$$
	\esp{d_{\mathcal{W}^{\delta}(1)}\pa{\hat{f}_J,f}} \leq C_{A,d,\delta,\beta,R} \pa{n\alpha^2}^{\frac{-(\beta+\delta)}{2\beta+2d}}.
	$$
	In this case, we recover the optimal rate, as it coincides with the lower bound obtained for all $\delta<d$. 
	
	\item[Case 2:] 
	If $\delta>d/2$, then the sum converges as $J$ goes to infinity, and we obtain that 
	$$
	\tilde{\Sigma}_J \leq C_{d,\delta}\frac{J^{d/2}}{\sqrt{n\alpha^2}}.
	$$
	Thus, the variance and bias terms are of the same order if
	$J^{-(\beta+\delta)} \asymp \frac{J^{d/2}}{\sqrt{n\alpha^2}},$
	that is if 
	$J \asymp \pa{n\alpha^2}^{\frac{1}{2\beta+2\delta+d}}.$
	In that case,  we obtain
	$$
	\esp{d_{\mathcal{W}^{\delta}(1)}\pa{\hat{f}_J,f}} \leq C_{A,d,\delta,\beta,R} \pa{n\alpha^2}^{\frac{-(\beta+\delta)}{2\beta+2\delta+d}}.
	$$
	In this case, for all $d/2<\delta<d$, the upper bound is larger than the lower bound obtained for all $\delta<d$ (which is the rate in Case 1). 
	Moreover, as $d>0$, for all $\delta>d$, the upper bound is also larger than $(n\alpha^2)^{-1/2}$ which is the minimax rate in that case. 
	Hence, when $\delta>d/2$, we never recover the optimal rate. 
	This comes from the fact that, when we do not use blocks, the variance of the private data is too high as it is proportional to the number of privatized coefficients. 
	
	\item[Case 3:] If $\delta = d/2$, a comparison between sum and integral leads to 
	$$\sum_{j=1}^J \frac{1}{j} \leq 1+\int_{1}^J \frac{dx}{x} = 1 + \log(J),$$
	and thus
	$$\tilde\Sigma_J \leq C_{d,\delta}\frac{(J\log(J))^{d/2}}{\sqrt{n\alpha^2}}.$$
	Then setting 
	$$
	J \asymp \pa{\frac{n\alpha^2}{\log(n\alpha^2)^{d}}}^{\frac{1}{2\beta+2\delta+d}},
	$$
	leads to 
	$$
	\tilde\Sigma_J 
	\leq C_{d,\delta,\beta}\frac{(n\alpha^2)^{\frac{d/2}{2\beta+2\delta+d}-1/2}}{\log(n\alpha^2)^{\frac{d}{2\beta+2\delta+d}\times\frac{d}{2}}} \log(J)^{d/2}
	\leq C_{d,\delta,\beta}\pa{\frac{n\alpha^2}{\log(n\alpha^2)^d}}^{\frac{-(\beta+\delta)}{2\beta+2\delta+d}},
	$$
	since 
	$
	\log(J) \asymp \frac{1}{2\beta+2\delta+d}\cro{\log(n\alpha^2) - d\log(n\alpha^2)} 
	\leq C_{d,\delta,\beta}\log(n\alpha^2).
	$
	In that case, we obtain that the bias and the variance are of the same order, and
	$$
	\esp{d_{\mathcal{W}^{\delta}(1)}\pa{\hat{f}_J,f}} \leq C_{A,d,\delta,\beta,R} \pa{\frac{n\alpha^2}{\log(n\alpha^2)^{d}}}^{\frac{-(\beta+\delta)}{2\beta+2\delta+d}}.
	$$
	As in the Case 2, this rate is not optimal, since it does not match with the lower bound for $\delta<d$ (not even up to a logarithmic term). 
\end{description}

\newpage
\bibliographystyle{apalike}

\nocite{*}
\bibliography{PrivateDensityEstim}

\end{document}